\documentclass[10pt]{amsart}
\usepackage{amscd}
\usepackage{amssymb}
\usepackage{latexsym}
\usepackage{graphicx}
\newcommand{\fourierabs}[1]{\lfloor #1 \rfloor}
\newcommand{\vangle}{\measuredangle}
\newcommand{\iterate}[1]{^{(#1)}}
\newcommand{\stFT}{\,\, \widetilde{} \,\,}
\newcommand{\FT}{\,\, \widehat{} \,\,}
\newcommand{\init}{\vert_{t = 0}}

\newcommand{\abs}[1]{\left\vert #1 \right\vert}
\newcommand{\fixedabs}[1]{\vert #1 \vert}
\newcommand{\bigabs}[1]{\bigl\vert #1 \bigr\vert}

\newcommand{\norm}[1]{\left\Vert #1 \right\Vert}
\newcommand{\fixednorm}[1]{\Vert #1 \Vert}
\newcommand{\bignorm}[1]{\bigl\Vert #1 \bigr\Vert}

\newcommand{\twonorm}[2]{\norm{#1}_{L^2#2}}

\newcommand{\C}{\mathbb{C}}

\newcommand{\R}{\mathbb{R}}
\newcommand{\Z}{\mathbb{Z}}

\newcommand{\innerprod}[2]{\left\langle \, #1 , #2 \, \right\rangle}

\newcommand{\angles}[1]{\langle #1 \rangle}

\DeclareMathOperator{\diag}{diag}

\DeclareMathOperator{\im}{Im}

\newtheorem{theorem}{Theorem}

\newtheorem{lemma}{Lemma}
\newtheorem{corollary}{Corollary}

\newtheorem*{theorem1}{Theorem 1}

\theoremstyle{definition}

\theoremstyle{remark}
\newtheorem{remark}{Remark}

\title[DKG in two space dimensions]{Local well-posedness below the charge norm for the Dirac-Klein-Gordon system in two space dimensions}

\author[P. D'Ancona]{Piero D'Ancona}

\address{Department of Mathematics\\
University of Rome ``La Sapienza''\\
Piazzale Aldo Moro 2\\
I-00185 Rome\\ Italy}
\email{dancona@mat.uniroma1.it}

\author[D. Foschi]{Damiano Foschi}
\address{Department of Pure and Applied Mathematics\\
University of L'Aquila\\
Via Vetoio, loc. Coppito\\
I-67010 L'Aquila\\ Italy}
\email{foschi@univaq.it}

\author[S. Selberg]{Sigmund Selberg}
\address{Department of Mathematical Sciences\\
Norwegian University of Science and Technology\\
Alfred Getz' vei 1\\
N-7491 Trondheim\\ Norway}
\thanks{The last author was supported by the Research Council of Norway, project no.\ 160192/V30, PDE and Harmonic Analysis.}
\email{sigmund.selberg@math.ntnu.no}


\begin{document}

\begin{abstract} 
We prove that the Cauchy problem for the Dirac-Klein-Gordon equations in two space dimensions is locally well-posed in a range of Sobolev spaces of negative index for the Dirac spinor, and an associated range of spaces of positive index for the meson field. In particular, we can go below the charge norm, that is, the $L^2$ norm of the spinor. We hope that this can have implications for the global existence problem, since the charge is conserved. Our result relies on the null structure of the system, and bilinear space-time estimates for the homogeneous wave equation.
\end{abstract}

\maketitle


\section{Introduction}\label{Section1}

We study the coupled Dirac-Klein-Gordon system of equations (DKG), which reads
\begin{equation}\label{DKG1}
\left\{
\begin{alignedat}{2}
  &\left( -i \gamma^\mu \partial_\mu + M \right) \psi = \phi \psi &\qquad& \left( M \ge 0 \right),
  \\
  &\left( - \square + m^2 \right) \phi = \psi^\dagger \gamma^0 \psi  &\qquad& \left( \square = - \partial_t^2 + \Delta, \,\, m \ge 0 \right),
\end{alignedat}
\right.
\end{equation}
where $\phi : \R^{1+n} \to \R$ represents a meson field and $\psi : \R^{1+n} \to \C^N$ is the Dirac spinor field, regarded as a column vector in $\C^N$; the dimension $N$ of the spin space depends on the space dimension $n$. Points in the Minkowski space-time $\R^{1+n}$ are written $(t,x)$, where $x = (x^1,\dots,x^n)$; we also denote $t=x^0$ when convenient. We write $\partial_\mu = \partial/\partial x^\mu$ and $\partial_t = \partial_0$. Greek indices $\mu,\nu,\dots$ range over $0,1,\dots,n$, roman indices $j,k,\dots$ over $1,\dots,n$, and repeated upper and lower indices are implicitly summed over these ranges. The $\gamma^\mu$'s are $N \times N$ matrices which should satisfy
$$
  \gamma^\mu \gamma^\nu + \gamma^\nu\gamma^\mu = 2g^{\mu\nu} I,
  \qquad
  (\gamma^0)^\dagger = \gamma^0
  \qquad \text{and}
  \qquad (\gamma^j)^\dagger = - \gamma^j,
$$
where $g^{\mu\nu} = \diag(1,-1,\dots,-1)$. The superscript $\dagger$ denotes conjugate transpose.

We shall refer to the cases $n = 1$, $2$ and $3$ as 1d, 2d and 3d, respectively. In 3d, the smallest possible dimension of the spin space, i.e., the smallest $N$ for which a realization of the Dirac matrices $\gamma^\mu$ can be found, is $N = 4$, whereas in 2d and 1d, $N = 2$. In this article, we are interested in the 2d case.

Concerning the Cauchy problem, global existence in 1d was established by Chadam \cite{Chadam:1973} (see also~\cite{Bournaveas:2000,Fang:2004}), but remains open in space dimension two and higher. Motivated by the work of Klainerman and Machedon~\cite{Klainerman:1994b,Klainerman:1995a} on the Maxwell-Klein-Gordon and Yang-Mills equations, we want to attack the global existence problem by improving the local well-posedness theory and then make use of conserved quantities, but the snag is that the energy density of DKG does not have a definite sign (see~\cite{Glassey:1979}). However, a partial replacement may be the charge conservation:
$$
  \int \abs{\psi(t,x)}^2 \, dx = \mathrm{const.},
$$
which could prove useful for the global problem, provided one has local well-posedness with $L^2$ data for the Dirac spinor. In fact, Bournaveas~\cite{Bournaveas:2000} used the charge conservation to give a new proof of Chadam's global result in 1d, by proving a low regularity local well-posedness theorem. 

These considerations motivate our interest in the local well-posedness of the Cauchy problem for DKG with initial data
\begin{equation}\label{DKGdata}
  \psi \init = \psi_0 \in H^s, \qquad \phi \init = \phi_0 \in H^r, \qquad \partial_t \phi \init = \phi_1 \in H^{r-1},
\end{equation}
for minimal $s,r \in \R$. Here $H^s = H^s(\R^n)$ is the standard $L^2$-based Sobolev space. The corresponding homogeneous space will be denoted $\dot H^s$. To get an idea of the minimal regularity required for local well-posedness, note that in the massless case $M = m = 0$, DKG is invariant under the rescaling
$$
  \psi(t,x) \longrightarrow \lambda^{3/2} \psi(\lambda t,\lambda x), \qquad
  \phi(t,x) \longrightarrow \lambda \phi(\lambda t,\lambda x).
$$
The scale invariant data space is therefore
$$
  (\psi_0,\phi_0,\phi_1) \in \dot H^{(n-3)/2} \times \dot H^{(n-2)/2} \times \dot H^{(n-4)/2},
$$
and one does not expect well-posedness with any less regularity than this. Since the charge corresponds to the $L^2$ norm of the spinor, we may say that DKG is charge-critical in 3d and charge-subcritical in 2d and 1d.

On the other hand, DKG is a system of nonlinear wave equations with quadratic nonlinearities, and it is well-known (see~\cite{Lindblad:1996,Ponce:1993}) that for such equations there is in general a gap between the regularity predicted by scaling and the minimal regularity at which one has local well-posedness, and this gap increases as the space dimension decreases. This is due to buildup effects in the product terms, but if the quadratic nonlinearities satisfy Klainerman's null condition, the worst interactions of products of waves are cancelled, and less regularity is required for local well-posedness. This idea first appeared in~\cite{Klainerman:1993}.

For DKG, the complete null structure was established by the authors in ~\cite{Selberg:2006b}, building on earlier work by Klainerman and Machedon~\cite{Klainerman:1994a} and Beals and B\'ezard~\cite{Beals:1996}. The new idea developed in~\cite{Selberg:2006b} is that the quadratic form $\psi^\dagger \gamma^0 \psi$, which appears in the Klein-Gordon part of DKG, and which was already known to be a null form (see ~\cite{Klainerman:1994a,Beals:1996}), appears also in the Dirac part of DKG, not directly but via a duality argument.

In~\cite{Selberg:2006b} we used the null structure, combined with bilinear space-time estimates of Klainerman-Machedon type, to prove local well-posedness of DKG in 3d for data \eqref{DKGdata} with $s > 0$ and $r = s+1/2$. Thus, we get arbitrarily close to the scaling regularity. For earlier work on the 3d problem, see~\cite{Bournaveas:1999,Fang:2005}.

\begin{figure}[h]
   \centering
   \includegraphics{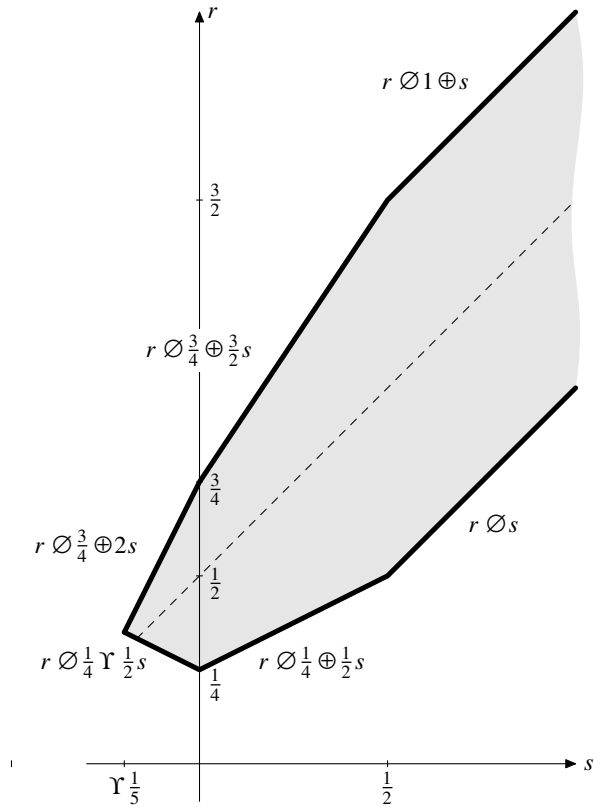}
   \caption{Local well-posedness holds in the interior of the shaded region, extending to the right. The dashed line $r=1/2+s$ represents the regularity predicted by scaling.}
   \label{fig:1}
\end{figure}

In the present work we study the 2d problem. Here it is harder to get close to the scaling regularity, since the range of space-time estimates narrows as the dimension decreases. On the other hand, DKG is charge-subcritical in 2d, so we stand a better chance of exploiting the charge conservation than in the charge-critical 3d case. In ~\cite{Bournaveas:2001}, Bournaveas proved local well-posedness of DKG in 2d for data \eqref{DKGdata} with $s>1/4$ and $r = s+1/2$ by using linear Strichartz type estimates (see also~\cite{Ponce:1993}), and also with $s>1/8$ and $r=s+5/8$ by using the null structure reported in~\cite{Klainerman:1994a}. Here we prove the following:

\begin{theorem}\label{Thm1} Suppose $(s,r) \in \R^2$ belongs to the convex region described by (see Figure \ref{fig:1})
$$
  s > -\frac{1}{5},
  \qquad
  \max\left(\frac{1}{4}-\frac{s}{2},\frac{1}{4}+\frac{s}{2},s\right)
  < r <
  \min\left(\frac{3}{4}+2s,\frac{3}{4}+\frac{3s}{2},1+s\right).
$$
Then the DKG system in 2d is locally well-posed for data \eqref{DKGdata}.
\end{theorem}

See Section \ref{NullBilinear} for a more precise statement. The proof relies on the null structure of the system, and some bilinear space-time estimates for the wave equation, to set up a contraction in $X^{s,b}$ type spaces.

The critical regularity for the spinor is $s = -1/2$, so one may ask what happens in the interval $-1/2 \le s \le -1/5$, which is not covered by our theorem. Partial answers to this question can be obtained by studying in more detail the iterates of the problem. We have some results for the first and second nontrivial iterates, which which will appear in a forthcoming paper. These results suggest, in particular, that there is no well-posedness for $s < -1/4$. This gap phenomenon is typical of 2d (and 1d) problems.

In connection with the iterates, we remark that the regularity of the first iterate for $\phi$ was studied by Zheng \cite{Zheng:1993}, who proved that if $\psi_0 \in L^2$, and $\psi\iterate{0}$ is the homogeneous part of $\psi$, then the modified iterate, defined by
$$
  \left(-\square+m^2\right){\Phi\iterate{1}} = g\innerprod{\beta \psi\iterate{0}}{\psi\iterate{0}},
  \qquad {\Phi\iterate{1}}(0,x) = \partial_t{\Phi\iterate{1}}(0,x) = 0,
$$
satisfies ${\Phi\iterate{1}}(t) \in L^2$ for $t > 0$, provided $g=g(t)$ is $C^1$ and $g(0) = 0$, the point being that energy and Sobolev estimates are enough to show this if $\psi_0 \in H^\varepsilon$ for some $\varepsilon > 0$, but they fail to give the result for $\varepsilon = 0$. However, our result here shows that the regularity is in fact far better than $L^2$, namely ${\Phi\iterate{1}}(t) \in H^{3/4-\varepsilon}$ for all $\varepsilon > 0$, if $\psi_0 \in L^2$; this follows from Theorem \ref{Thm2} below.

Throughout the rest of this article, the space dimension is understood to be $n = 2$. As a matter of convenience, we consider only the massless case $M=m=0$, but the discussion can easily be modified to include the linear mass terms, since we deal with local-in-time theory in a contraction mapping setup.

For convenience we rewrite the system \eqref{DKG1} in the form
\begin{equation}\label{DKG}
\left\{
\begin{aligned}
  &i \bigl( \partial_t + \alpha\cdot \nabla \bigr) \psi = - \phi \beta \psi,
  \\
  &\square \phi = - \innerprod{\beta \psi}{\psi},
\end{aligned}
\right.
\end{equation}
where $\alpha^j = \gamma^0\gamma^j$, $\beta = \gamma^0$, $\psi= [\psi_1,\psi_2]^T$, $\innerprod{z}{w} = w^\dagger z$ for column vectors $z,w \in \C^2$, $\nabla = (\partial_1,\partial_2)$ and $\alpha \cdot \nabla = \alpha^1\partial_1+\alpha^2\partial_2$. The Dirac matrices $\alpha^j, \beta$ should satisfy
\begin{equation}
  \label{DiracIdentity1}
  \beta^\dagger=\beta, \quad (\alpha^j)^\dagger = \alpha^j,
  \quad \beta^2 = (\alpha^j)^2 = I, \quad \alpha^j \beta + \beta \alpha^j = 0,
  \quad \alpha^j \alpha^k + \alpha^k \alpha^j = 2\delta^{jk} I.
\end{equation}
A particular representation in 2d is
$$
  \alpha^{1} = \sigma^{1} =
  \begin{pmatrix}
    0 & 1  \\
    1 & 0
  \end{pmatrix},
  \qquad
  \alpha^2 = \sigma^{2} =
  \begin{pmatrix}
    0 & -i  \\
    i & 0
  \end{pmatrix},
  \qquad
  \beta = \sigma^{3} =
  \begin{pmatrix}
    1 & 0  \\
    0 & -1
  \end{pmatrix},
$$
where the $\sigma$'s are the Pauli matrices. Thus, $\psi^\dagger \beta \psi = \innerprod{\beta \psi}{\psi} = \fixedabs{\psi_1}^2 - \fixedabs{\psi_2}^2$.

The rest of this paper is organized as follows: In Section \ref{Notation} we fix the notation, and define the function spaces in which we iterate. We also review the splitting of the Dirac equation into positive and negative energy parts. In Section \ref{NullBilinear}, we reduce Theorem \ref{Thm1} to proving two bilinear $X^{s,b}$ estimates, we review the crucial null structure of the bilinear forms involved, and we discuss two key ingredients that will be used to prove the $X^{s,b}$ estimates: a bilinear generalization of the Strichartz estimate for free waves in 2d, and a variation on the Klainerman-Machedon estimate in 2d. In Sections \ref{Thm2Proof} and \ref{Thm3Proof} we prove the $X^{s,b}$ estimates, in Section \ref{Optimality} we prove that they are optimal up to endpoint cases, and in Section \ref{Thm4Proof} we prove the bilinear Strichartz type estimate.

\section{Notation}\label{Notation}

In estimates, we use $C$ to denote a large, positive constant which can change from line to line. If $C$ is absolute, or only depends on parameters which are considered fixed, then we often write $\lesssim$, which means $\le$ up to multiplication by $C$. If $X,Y$ are nonnegative quantities, $X \sim Y$ means $C^{-1}Y \le X \le CY$ for some absolute constant $C \gg 1$. Naturally, we then define $X \ll Y$ to mean $X \le C^{-1} Y$, and $X \gg Y$ to mean $X \ge CY$. Throughout we use the notation $\angles{\cdot} = 1 + \abs{\cdot}$. The characteristic function of a set $A$ is denoted $\chi_A$.

All $L^p$ norms are with respect to Lebesgue measure on $\R^2$, $\R^{1+2}$ or $\R$, unless stated otherwise. Often we indicate by a subscript which variable or variables the norm is taken over, as in $L^p_x$, $L^p_{t,x}$ or $L^p_t$. The norm on $L^q_t\bigl(\R;L^r_x(\R)\bigr)$ is denoted
$$
  \norm{u}_{L_t^q L_x^r} = \norm{ \norm{u(t,x)}_{L_x^r} }_{L_t^q}
  = \left( \int_{-\infty}^\infty \left( \int_{\R^2} \abs{u(t,x)}^r \, dx \right)^{q/r} \, dt \right)^{1/q},
$$
with the usual modification if $q$ or $r$ equals $\infty$.

The Fourier transforms in space and space-time are defined by
$$
  \widehat f(\xi) = \int_{\R^2} e^{-ix\cdot\xi} f(x) \, dx,
  \qquad
  \widetilde u(\tau,\xi) = \int_{\R^{1+2}} e^{-i(t\tau+x\cdot\xi)} u(t,x) \, dt \, dx.
$$
Then Plancherel's theorem comes out somewhat awkwardly as $\twonorm{\widehat f}{} = 2\pi \twonorm{f}{}$ and $\twonorm{\widetilde u}{} = (2\pi)^{3/2} \twonorm{u}{}$. To avoid having to keep track of irrelevant normalization factors, we use $\simeq$ to mean equality up to multiplication by some fixed, positive factor. Given $u(t,x)$, we denote by $\fourierabs{u}$ the function whose space-time Fourier transform is $\abs{\widetilde u}$.

We write $D = \nabla/i$, where $\nabla = (\partial_1,\partial_2)$. Then $(Df)\FT(\xi) = \xi \widehat f(\xi)$, which explains the notation $\phi(D)$ for the multiplier defined by
$$
  \widehat{\phi(D) f}(\xi) = \phi(\xi) \widehat f(\xi),
$$
for a given symbol $\phi$. The multipliers $\abs{D}^s$ and $\angles{D}^s = (1+\abs{D})^s$ are used to define $\dot H^s$ and $H^s$ as the completions of the Schwartz space $\mathcal S(\R^2)$ with respect to the norms
$$
  \norm{f}_{\dot H^s} = \norm{\abs{D}^s f}_{L^2_x} \simeq \norm{\abs{\xi}^s f}_{L^2_\xi},
  \qquad
  \norm{f}_{H^s} = \norm{\angles{D}^s f}_{L^2_x} \simeq \norm{\angles{\xi}^s f}_{L^2_\xi}.
$$
Note that for $\dot H^s$ we need $s > -1$, since the space dimension is $n=2$.

The operator $-i\alpha \cdot \nabla$ appearing in the Dirac equation is awkward to deal with, since it mixes the components of the spinor it acts on. To simplify, we decompose the spinor along an eigenbasis of the operator. Specifically, the matrix-valued symbol of $-i\alpha \cdot \nabla$ is $\xi \cdot \alpha = \xi_1 \alpha^1 + \xi_2 \alpha^2$, which is hermitian and satisfies $(\xi\cdot\alpha)^2 = \abs{\xi}^2 I$, on account of \eqref{DiracIdentity1}. Thus, the eigenvalues are $\pm \abs{\xi}$, and the corresponding projections onto the one-dimensional eigenspaces are
\begin{equation}\label{Projection}
  \Pi_{\pm}(\xi) = \frac{1}{2} \left( I \pm \frac{\xi}{\abs{\xi}} \cdot \alpha \right)
  =
  \frac{1}{2}
  \begin{pmatrix}
        1 & \pm(\hat\xi_1-i\hat\xi_2) \\
        \pm(\hat\xi_1+i\hat\xi_2) & 1 \\
     \end{pmatrix},
  \qquad \text{where} \quad \hat\xi \equiv \frac{\xi}{\abs{\xi}}.
\end{equation}
Then  $-i\alpha \cdot \nabla = \abs{D} \Pi_+(D) - \abs{D} \Pi_-(D)$, so the solution of the linear Cauchy problem
\begin{equation}\label{LinearDiracIVP}
  -i\bigl( \partial_t + \alpha\cdot \nabla \bigr) \psi = F, \qquad \psi(0,x) = \psi_0(x),
\end{equation}
splits into $\psi = \psi_+ + \psi_-$, where $\psi_\pm = \Pi_\pm(D) \psi$ satisfy
\begin{equation}\label{LinearDiracIVPSplit}
\left\{
\begin{alignedat}{2}
  &\bigl(-i\partial_t \pm \abs{D}\bigr) \psi_\pm = F_\pm
  &\qquad &\left( \psi_\pm = \Pi_\pm(D) \psi, \,\, F_\pm = \Pi_\pm(D) F \right),
  \\
  &\psi_\pm(0,x) = \psi_0^\pm(x)
  & &\left( \psi_0^\pm = \Pi_\pm(D) \psi_0 \right).
\end{alignedat}
\right.
\end{equation}
Note that in a physical interpretation, at least for the free case $F = 0$, the spinors $\psi_+$ and $\psi_-$ correspond to positive and negative energies, respectively. The free propagator for $-i\partial_t \pm \abs{D}$ is the multiplier $S_\pm(t) = e^{\mp it\abs{D}}$ with symbol $e^{\mp it \abs{\xi}}$. Note that $S_\pm(t)$ acts componentwise on spinors.

To prove Theorem \ref{Thm1} we shall iterate $\psi_\pm$ and $\phi$ in $X^{s,b}$ type spaces associated to the operators $-i\partial_t \pm \abs{D}$ and $\square$, whose symbols are $\tau \pm \abs{\xi}$ and $\tau^2-\abs{\xi}^2$, respectively. See \cite{Selberg:2006b} for more details about the following spaces. Let $D_\pm$ be the multipliers with symbols $\abs{\tau}\pm\abs{\xi}$. For $s,b \in \R$, we define $X^{s,b}_\pm$, $H^{s,b}$ and $\mathcal H^{s,b}$ to be the completions of the Schwartz space $\mathcal S(\R^{1+2})$ with respect to the norms
\begin{align*}
  \norm{u}_{X^{s,b}_\pm} &= \bignorm{\angles{D}^s \angles{-i\partial_t \pm \abs{D}}^b u}_{L^2_{t,x}}
  \simeq \bignorm{\angles{\xi}^s \angles{\tau\pm\abs{\xi}}^b \widetilde u(\tau,\xi)}_{L^2_{\tau,\xi}},
  \\
  \norm{u}_{H^{s,b}} &= \bignorm{\angles{D}^s \angles{D_-}^b u}_{L^2_{t,x}}
  \simeq \bignorm{\angles{\xi}^s \angles{\abs{\tau} - \abs{\xi}}^b \widetilde u(\tau,\xi)}_{L^2_{(\tau,\xi)}},
  \\
  \norm{u}_{\mathcal H^{s,b}} &= \norm{u}_{H^{s,b}} + \norm{\partial_t u}_{H^{s-1,b}}
  \sim \bignorm{\angles{D}^{s-1} \angles{D_+} \angles{D_-}^b u}_{L^2_{t,x}},
\end{align*}
where $\angles{\cdot} = 1 + \abs{\cdot}$. We also need the restrictions to a time slab
$
  S_T = (0,T) \times \R^2,
$
since we study local solutions. The restriction $X_\pm^{s,b}(S_T)$ is a Banach space with norm
$$
  \norm{u}_{X_\pm^{s,b}(S_T)} = \inf \left\{
  \norm{v}_{X_\pm^{s,b}} : \text{$v \in X_\pm^{s,b}$ and $v = u$ on $S_T$}
  \right\}.
$$
Completeness holds since $X_\pm^{s,b}(S_T) = X_\pm^{s,b}/\mathcal M_\pm$, where $\mathcal M_\pm = \{ v \in X_\pm^{s,b} : \text{$v = 0$ on $S_T$} \}$ is closed. The restrictions $H^{s,b}(S_T)$ and $\mathcal H^{s,b}(S_T)$ are defined in the same way.
\section{Null structure and bilinear estimates}\label{NullBilinear}

The complete null structure of DKG, found recently in \cite{Selberg:2006b}, rests on the cancellation properties of the matrix-valued symbol
$$
  \sigma_{\pm_1,\pm_2}(\eta,\zeta) = \Pi_{\pm_2}(\zeta)\beta\Pi_{\pm_1}(\eta)
  = \beta\Pi_{\mp_2}(\zeta)\Pi_{\pm_1}(\eta),
$$
where to get the last equality we used $\beta \Pi_\pm(\xi) = \Pi_{\mp}(\xi) \beta$, which follows from \eqref{DiracIdentity1}. By orthogonality, $\Pi_{\mp_2}(\zeta)\Pi_{\pm_1}(\eta)$ vanishes when the vectors $\pm_1\eta$ and $\pm_2\zeta$ line up in the same direction. The following lemma, proved in \cite{Selberg:2006b}, quantifies this cancellation. We shall use the notation $\vangle(\eta,\zeta)$ for the angle between vectors $\eta,\zeta \in \R^2$.

\begin{lemma}\label{NullLemma}
$\sigma_{\pm_1,\pm_2}(\eta,\zeta) = O \left( \vangle(\pm_1\eta,\pm_2\zeta) \right)$.
\end{lemma}

Note the following consequence: If $\psi,\psi' : \R^{1+2} \to \C^2$ are Schwartz functions, then
\begin{equation}\label{NullFormEstimate}
\begin{aligned}
  &\abs{\innerprod{\beta \Pi_{\pm_1}(D) \psi}{\Pi_{\pm_2}(D) \psi'}\stFT(\tau,\xi)}
  \\
  &\qquad\qquad\le
  \int_{\R^{1+2}} \abs{\innerprod{\beta\Pi_{\pm_1}(\eta)\widetilde \psi(\lambda,\eta)}
  {\Pi_{\pm_2}(\eta-\xi)\widetilde \psi'(\lambda-\tau,\eta-\xi)}} \, d\lambda \, d\eta
  \\
  &\qquad\qquad\lesssim \int_{\R^{1+2}} \theta_{\pm_1,\pm_2} \abs{\widetilde \psi(\lambda,\eta)}
  \abs{\widetilde \psi'(\lambda-\tau,\eta-\xi)} \, d\lambda \, d\eta,
\end{aligned}
\end{equation}
where $\theta_{\pm_1,\pm_2} = \vangle\bigl(\pm_1\eta,\pm_2(\eta-\xi)\bigr)$. Here we used the self-adjointness of the projections to move $\Pi_{\pm_2}(\eta-\xi)$ in front of $\beta$, and then we applied the lemma.

Let us now restate our main theorem in a more precise form.

\begin{theorem1}
Suppose $(s,r) \in \R^2$ belongs to the convex region described by (see Figure \ref{fig:1})
$$
  s > -\frac{1}{5},
  \qquad
  \max\left(\frac{1}{4}-\frac{s}{2},\frac{1}{4}+\frac{s}{2},s\right)
  < r <
  \min\left(\frac{3}{4}+2s,\frac{3}{4}+\frac{3s}{2},1+s\right).
$$
Then for any data $\psi \init = \psi_0 \in H^s$, $\phi \init = \phi_0 \in H^r$ and $\partial_t \phi \init = \phi_1 \in H^{r-1}$,
there exist a time $T > 0$, depending continuously on the $H^s \times H^r \times H^{r-1}$ norm of the data, and a solution
$$
  \psi \in C([0,T];H^s), \qquad \phi \in C([0,T];H^r) \cap C^1([0,T];H^{r-1}),
$$
of the DKG system \eqref{DKG} on $(0,T) \times \R^2$, satisfying the initial condition above. Furthermore, writing $\psi = \psi_+ + \psi_-$, where $\psi_\pm = \Pi_\pm(D) \psi$, the solution has the regularity
$$
  \psi_\pm \in X_\pm^{s,1/2+\varepsilon}(S_T),
  \qquad \phi \in \mathcal H^{r,1/2+\varepsilon}(S_T)
  \qquad \text{for $\varepsilon = \varepsilon(s,r) > 0$ sufficiently small.}
$$
Moreover, the solution is unique in this class, and depends continuously on the data.
\end{theorem1}

The first step in the proof is to use \eqref{LinearDiracIVP} and \eqref{LinearDiracIVPSplit} to reformulate the system \eqref{DKG} as
$$
\left\{
\begin{aligned}
  &\bigl( -i\partial_t \pm \abs{D} \bigr) \psi_\pm = - \Pi_\pm(D) \left( \phi \beta \psi \right),
  \\
  &\square \phi = - \innerprod{\beta \psi}{\psi}.
\end{aligned}
\right.
$$
Then by a standard iteration argument, using estimates for the Cauchy problem \eqref{LinearDiracIVPSplit} in $X^{s,1/2+\varepsilon}_\pm(S_T)$, and for the Cauchy problem for the wave equation in $\mathcal H^{r,1/2+\varepsilon}(S_T)$, see \cite[Lemmas 5 and 6]{Selberg:2006b}, the theorem reduces to the following bilinear estimates:
\begin{align}
  \label{BilinearA}
  \norm{\Pi_{\pm_2}(D) \left(\phi\beta\Pi_{\pm_1}(D) \psi\right)}_{X_{\pm_2}^{s,-1/2+2\varepsilon}(S_T)}
  &\lesssim \norm{\phi}_{H^{r,1/2+\varepsilon}(S_T)} \norm{\psi}_{X_{\pm_1}^{s,1/2+\varepsilon}(S_T)},
  \\
  \label{BilinearB}
  \norm{\innerprod{\beta\Pi_{\pm_1}(D) \psi}{\Pi_{\pm_2}(D)\psi'}}_{H^{r-1,-1/2+2\varepsilon}(S_T)}
  &\lesssim \norm{\psi}_{{X_{\pm_1}^{s,1/2+\varepsilon}(S_T)}}
  \norm{\psi'}_{X_{\pm_2}^{s,1/2+\varepsilon}(S_T)},
\end{align}
for all combinations of signs, where $S_T = (0,T) \times \R^2$, and we assume $0 < T \le 1$. Here and in the rest of the paper, it is understood that implicit constants may depend on $s$,$r$ and $\varepsilon$, but not on $T$.

We shall prove the following theorems, optimal up to endpoint cases, which imply Theorem \ref{Thm1} (see Remark \ref{CutoffRemark} below). Here it is understood that $\phi$ is real-valued and $\psi,\psi'$ are $\C^2$-valued.

\begin{theorem}\label{Thm2}
Suppose $s > -1/4$ and $r < \min(3/4+2s,3/4+3s/2,1+s)$. Then there exists $\varepsilon > 0$ such that \begin{equation}\label{BilinearBB}
  \norm{\innerprod{\beta\Pi_{\pm_1}(D) \psi}{\Pi_{\pm_2}(D)\psi'}}_{H^{r-1,-1/2+2\varepsilon}}
  \lesssim \norm{\psi}_{{X_{\pm_1}^{s,1/2+\varepsilon}}}
  \norm{\psi'}_{X_{\pm_2}^{s,1/2+\varepsilon}}
\end{equation}
for all $\psi,\psi' \in \mathcal S(\R^{1+2})$ such that $\psi$ is supported in $[-2,2] \times \R^2$. Moreover, the estimate fails if $s < -1/4$ or $r > \min(3/4+2s,3/4+3s/2,1+s)$.
\end{theorem}

\begin{theorem}\label{Thm3}
Suppose $s \in \R$ and $r > \max(-s,1/4-s/2,1/4+s/2,s)$. Then there exists $\varepsilon > 0$ such that
\begin{equation}\label{BilinearAA}
  \norm{\Pi_{\pm_2}(D) \left(\phi\beta\Pi_{\pm_1}(D) \psi\right)}_{X_{\pm_2}^{s,-1/2+2\varepsilon}}
  \lesssim \norm{\phi}_{H^{r,1/2+\varepsilon}} \norm{\psi}_{X_{\pm_1}^{s,1/2+\varepsilon}}
\end{equation}
for all $\phi,\psi \in \mathcal S(\R^{1+2})$ such that $\psi$ is supported in $[-2,2] \times \R^2$. Moreover, the estimate fails if $r < \max(-s,1/4-s/2,1/4+s/2,s)$.
\end{theorem}

In \cite{Selberg:2006b} it was shown that \eqref{BilinearAA} is equivalent, by duality, to an estimate similar to \eqref{BilinearBB}, namely
\begin{equation}\label{BilinearAAA}\tag{$\text{\ref{BilinearAA}}'$}
  \norm{\innerprod{\beta\Pi_{\pm_1}(D) \psi}{\Pi_{\pm_2}(D)\psi'}}_{H^{-r,-1/2-\varepsilon}}
  \lesssim \norm{\psi}_{{X_{\pm_1}^{s,1/2+\varepsilon}}}
  \norm{\psi'}_{X_{\pm_2}^{-s,1/2-2\varepsilon}},
\end{equation}
which must hold for all $\psi, \psi' \in \mathcal S(\R^{1+2})$ such that $\psi$ is supported in $[-2,2] \times \R^2$. Note the advantage of this formulation, in that the null form appears again.


\begin{remark}\label{CutoffRemark} Theorem \ref{Thm2} implies \eqref{BilinearB}, as we now show. A similar argument can be used to show that Theorem \ref{Thm3} implies \eqref{BilinearA}. Fix a smooth cutoff $\chi(t)$ such that $\chi(t) = 1$ for $\abs{t} \le 1$ and $\chi(t) = 0$ for $\abs{t} \ge 2$. Let $\psi \in X_{\pm_1}^{s,1/2+\varepsilon}(S_T)$ and $\psi' \in X_{\pm_2}^{s,1/2+\varepsilon}(S_T)$. Viewing $X_\pm^{s,b}(S_T)$ as a space of equivalence classes, the equivalence relation being equality on $S_T$, we let $\Psi \in X_{\pm_1}^{s,1/2+\varepsilon}$ and $\Psi' \in X_{\pm_2}^{s,1/2+\varepsilon}$ denote arbitrary representatives of $\psi$ and $\psi'$, respectively. Recalling the assumption $0 < T \le 1$, we observe that Lemma \ref{Thm2} implies
\begin{align*}
  \norm{\innerprod{\beta\Pi_{\pm_1}(D) \psi}{\Pi_{\pm_2}(D)\psi'}}_{H^{r-1,-1/2+2\varepsilon}(S_T)}
  &\le
  \norm{\innerprod{\beta\Pi_{\pm_1}(D)(\chi\Psi)}{\Pi_{\pm_2}(D)\Psi'}}_{H^{r-1,-1/2+2\varepsilon}}
  \\
  &\lesssim \bignorm{\chi\Psi}_{{X_{\pm_1}^{s,1/2+\varepsilon}}}
  \bignorm{\Psi'}_{X_{\pm_2}^{s,1/2+\varepsilon}}
  \\
  &\lesssim \bignorm{\Psi}_{{X_{\pm_1}^{s,1/2+\varepsilon}}}
  \bignorm{\Psi'}_{X_{\pm_2}^{s,1/2+\varepsilon}},
\end{align*}
and taking the infimum over all representatives $\Psi,\Psi'$ yields \eqref{BilinearB}. In the last step we used the easily proved estimate
\begin{equation}\label{SimpleEstimate}
  \fixednorm{\chi u}_{X_\pm^{s,b}} \le C_{\chi,b} \norm{u}_{X_\pm^{s,b}},
\end{equation}
valid for $b \ge 0$.
\end{remark}

In addition to the null form estimate \eqref{NullFormEstimate}, the main tools for proving Theorems \ref{Thm2} and \ref{Thm3} are some bilinear spacetime estimates for 2d free waves, which we now discuss. Recall that $S_\pm(t) = e^{\mp it \abs{D}}$ is the free propagator for $-i\partial_t \pm \abs{D}$. In the following discussion we let $f,g \in \mathcal S(\R^2)$ and write
$$
  u(t) = u_\pm(t) = S_\pm(t)f,
  \qquad
  v(t) = v_\pm(t) = S_\pm(t)g.
$$
In estimates where the signs do not matter, we skip the subscript indicating the sign.

We begin with a generalization of the Strichartz estimate for free waves in 2d. The following extends (in the 2d case) an estimate due to Klainerman and Tataru \cite{Klainerman:1999}. 

\begin{theorem}\label{Thm4}
The estimate
$$
  \bignorm{\abs{D}^{-s_3} (uv)}_{L_t^q L_x^2} \lesssim \norm{f}_{\dot H^{s_1}} \norm{g}_{\dot H^{s_2}}
$$
holds if
$$
  \begin{cases}
  4 \le q \le \infty,
  \\
  s_1+s_2+s_3 = 1-1/q,
  \\
  s_1,s_2 < 1-1/q,
  \\
  s_1+s_2 > 1/q\quad( \iff s_3 < 1-2/q  ).
  \end{cases}
$$
\end{theorem}

The case where $q = 4$, $s_1=s_2$ and $s_3 \le 0$ was proved in \cite{Klainerman:1999}. The above theorem is sharp up to endpoint cases, in view of \cite[Proposition 14.15]{Foschi:2000}. In particular, there are no estimates of this form with $q < 4$.

Note that \eqref{BilinearAA} and \eqref{BilinearBB} are $L^2$ in both space and time, so to apply the above theorem we need to use the finite support in time. For example, if we restrict to a finite time interval $0 \le t \le T$, then by H\"older's inequality in time and Theorem \ref{Thm4} with $q=4$,
\begin{equation}\label{FiniteTime}
  \norm{ \abs{D}^{-s_3}(uv)}_{L^2([0,T] \times \R^2)} \lesssim T^{1/4} \norm{f}_{\dot H^{s_1}} \norm{g}_{\dot H^{s_2}}
  \qquad \text{if}
  \quad
  \begin{cases}
  s_1+s_2+s_3 = 3/4,
  \\
  s_1,s_2 < 3/4,
  \\
  s_1+s_2 > 1/4.
  \end{cases}
\end{equation}

Theorem \ref{Thm4} suffices to prove Theorem \ref{Thm3}, and also Theorem \ref{Thm2}, with the exception of one particularly delicate case (see subsection \ref{DifficultCase}), where we need to use a variant of the Klainerman-Machedon estimate, which we now discuss.

For $q < 4$ there are no estimates of the type considered in Theorem \ref{Thm4}, but Klainerman and Machedon \cite{Klainerman:1996} proved that if the product $uv$ is replaced by the null form $D_-^{\gamma}(uv)$, where $\gamma \ge 1/4$, and $D_-$ denotes the multiplier with symbol $\bigabs{ \abs{\tau} - \abs{\xi} }$, then one can obtain a range of estimates with $q=2$.

\begin{theorem}\label{KMthm} \cite{Klainerman:1996,Foschi:2000}. The estimate
$$
  \norm{ \abs{D}^{-s_3}D_-^{1/4}(uv)}_{L^2_{t,x}} \lesssim  \norm{f}_{\dot H^{s_1}} \norm{g}_{\dot H^{s_2}}
$$
holds if and only if
$$
  \begin{cases}
  s_1+s_2+s_3 = 3/4,
  \\
  s_1,s_2 < 3/4,
  \\
  s_1+s_2 > 1/2 \quad( \iff s_3 < 1/4  ).
  \end{cases}
$$
\end{theorem}

The fact that $\gamma \ge 1/4$ is required to have an $L^2$ space-time estimate for $D_-^{\gamma}(uv)$, is related to the gap phenomenon in 2d (there is a gap between the regularity predicted by scaling and the regularity needed to have local well-posedness). Indeed, the null forms in DKG correspond to $\gamma = 1/2$, and in principle the null symbol may then be completely cancelled against the weight of the $X^{s,-1/2}$ norm in which the null form is estimated, resulting in a corresponding gain of ``elliptic'' derivatives. But in reality we can only cancel ``half'' of the symbol, due to the restriction $\gamma \ge 1/4$ in the $L^2$ estimate. This corresponds to a gap of a quarter of a derivative down to the scaling regularity. This is exactly the loss we incur if instead we cancel the full null symbol and then, assuming finite support in time, apply H\"older's inequality in time to replace $L^2_t$ by $L^4_t$, at which point we are in a position to apply Theorem \ref{Thm4}. The latter approach is always better, however: comparing the last theorem with \eqref{FiniteTime}, we see that the conditions on the $s_i$ are the same, except that in \eqref{FiniteTime} one only needs $s_1+s_2 > 1/4$, as opposed to $s_1+s_2 > 1/2$ in Theorem \ref{KMthm}. 

As it stands, Theorem \ref{KMthm} is not useful in the present context, in view of the preceding remarks. However, as  observed in \cite[Theorem 6(b)]{Selberg:1999} (or see \cite[Theorem 12.1]{Foschi:2000}), the condition $s_1 + s_2 > 1/2$ in Theorem \ref{KMthm} can be relaxed to $s_1 + s_2 > 1/4$ for products of type $(+,+)$ and $(-,-)$, i.e., if $uv$ is $u_+v_+$ or $u_-v_-$. In fact,
\begin{equation}\label{ImprovedKMest}
  \norm{ \abs{D}^{-s_3}D_-^{1/4}(u_+ v_+)}_{L^2_{t,x}} \lesssim \norm{f}_{\dot H^{s_1}} \norm{g}_{\dot H^{s_2}}
   \qquad \text{if}
  \quad
  \begin{cases}
  s_1+s_2+s_3 = 3/4,
  \\
  s_1,s_2 < 3/4,
  \\
  s_1+s_2 > 1/4.
  \end{cases}
\end{equation}
We are not aware of any previous application of this improved estimate. Here we shall use it to get a decisive improvement of \eqref{FiniteTime}, in the case where the product is of type $(+,+)$ or $(-,-)$ with high frequencies interacting to give output at low frequency. To make this precise, we introduce a bilinear operator $(\cdot,\!\cdot)_{\text{HH$\to$L}}$ which isolates this interaction:
$$
  \left[(f,g)_{\text{HH$\to$L}}\right]\FT(\xi)
  = \int_{\R^2} \mathbb{1}_{\abs{\xi} \ll \fixedabs{\eta} + \fixedabs{\xi-\eta}}
  \widehat f(\eta) \widehat g(\xi-\eta)
  \, d\eta.
$$
Here $\mathbb{1}_{\abs{\xi} \ll \fixedabs{\eta} + \fixedabs{\xi-\eta}}$ is the characteristic function of the set $\{ \eta : \abs{\xi} \ll \fixedabs{\eta} + \fixedabs{\xi-\eta} \}$. Then for free waves $u_+(t) = S_+(t)f$ and $v_+(t) = S_+(t)g$, using the fact that $\widetilde u_+(\tau,\xi) = \delta(\tau+\abs{\xi})\widehat f(\xi)$ and $\widetilde v_+(\tau,\xi) = \delta(\tau+\abs{\xi})\widehat g(\xi)$, we have
\begin{align*}
  &\left[D_-^{1/4}(u_+,v_+)_{\text{HH$\to$L}}\right]\stFT(\tau,\xi)
  \\
  &\qquad= \int_{\R^{1+2}} \bigabs{\abs{\tau} - \abs{\xi}}^{1/4}
  \mathbb{1}_{\abs{\xi} \ll \fixedabs{\eta} + \fixedabs{\xi-\eta}}
  \widehat f(\eta) \widehat g(\xi-\eta)
  \delta(\lambda+\fixedabs{\eta})\delta(\tau-\lambda+\fixedabs{\xi-\eta})
  \, d\lambda \, d\eta
  \\
  &\qquad= \int_{\R^2} \left( \fixedabs{\eta} + \fixedabs{\xi-\eta} - \abs{\xi} \right)^{1/4} \mathbb{1}_{\abs{\xi} \ll \fixedabs{\eta} + \fixedabs{\xi-\eta}}
  \widehat f(\eta) \widehat g(\xi-\eta)
  \delta(\tau+\fixedabs{\eta}+\fixedabs{\xi-\eta})
  \, d\eta,
\end{align*}
hence
$$
  \norm{\abs{D}^{-s_3} (u_+,v_+)_{\text{HH$\to$L}}}_{L^2_{t,x}}
  \sim \norm{\abs{D}^{-s_3} D_-^{1/4} (\abs{D}^{-1/8} u_+,\abs{D}^{-1/8}v_+)_{\text{HH$\to$L}}}_{L^2_{t,x}},
$$
and from \eqref{ImprovedKMest} we then obtain:

\begin{theorem}\label{Thm6}
We have
$$
  \norm{\abs{D}^{-s_3} (u_+,v_+)_{\text{HH$\to$L}}}_{L^2_{t,x}} \lesssim \norm{f}_{\dot H^{s_1}} \norm{g}_{\dot H^{s_2}}
  \qquad
  \text{if}
  \quad
  \begin{cases}
  s_1+s_2+s_3 = 1/2,
  \\
  s_1,s_2 < 5/8,
  \\
  s_1+s_2 > 0.
  \end{cases}
$$
\end{theorem}

Note that this is a tremendous improvement over \eqref{FiniteTime}, for this particular interaction.
This estimate will be used to handle a particularly delicate case occurring in the proof of Theorem \ref{Thm2} (see subsection \ref{DifficultCase}), where Theorem \ref{Thm4} fails hopelessly.

To end this section, let us state the $X^{s,b}$ versions of the free wave estimates discussed above.
It is a general principle (see, e.g., \cite[Lemma 4]{Selberg:2006b}) that Strichartz type estimates for a free propagator imply corresponding estimates for the $X^{s,b}$ space associated to the propagator. Thus, Theorems \ref{Thm4} and \ref{Thm6} imply, respectively:

\begin{corollary}\label{BilinearCorollary}
Suppose $\varepsilon > 0$ and $q,s_1,s_2,s_3$ satisfy the hypotheses in Theorem \ref{Thm4}. Then
$$
  \bignorm{\abs{D}^{-s_3} (uv)}_{L_t^q L_x^2} \le C_{q,s_1,s_2,\varepsilon}
  \norm{\abs{D}^{s_1} u}_{X_\pm^{0,1/2+\varepsilon}}
  \norm{\abs{D}^{s_2} v}_{X_\pm^{0,1/2+\varepsilon}}
  \qquad \text{for $u,v \in \mathcal S(\R^{1+2})$},
$$
for all combinations of signs.
\end{corollary}

\begin{corollary}\label{BilinearCorollary2}
Suppose $\varepsilon > 0$ and $s_1,s_2,s_3$ satisfy the hypotheses in Theorem \ref{Thm6}. Then
$$
  \norm{\abs{D}^{-s_3} (u,v)_{\text{HH$\to$L}}}_{L^2_{t,x}} \lesssim \norm{u}_{X_+^{s_1,1/+\varepsilon}} \norm{v}_{X_+^{s_2,1/+\varepsilon}}
  \qquad \text{for $u,v \in \mathcal S(\R^{1+2})$}.
$$
Note carefully the equality of the signs in the norms on the right.
\end{corollary}

\section{Proof of Theorem \ref{Thm2}}\label{Thm2Proof}

Without loss of generality we take $\pm_1 = +$ and write $\pm_2 = \pm$. We assume $s > -1/2$ and
\begin{equation}\label{rConditions}
  r < \min\left(\frac{3}{2}+4s,\frac{3}{4}+2s,\frac{3}{4}+\frac{3s}{2},1+s\right).
\end{equation}
In fact, if $\pm$ is a minus, we need to assume $s > -1/4$, and then the condition $r < 3/2+4s$ is redundant. The parameter $\varepsilon > 0$ will be chosen sufficiently small, depending on $s$ and $r$. In the sequel, whenever we say that some condition involving $\varepsilon$ holds, we mean that it holds for all $\varepsilon > 0$ small enough.
We assume that $\psi,\psi' \in \mathcal S(\R^{1+2})$ and $\psi$ is supported in $[-2,2] \times \R^2$. We use the notation $\fourierabs{\psi}$ for the function whose Fourier transform is $\fixedabs{\widetilde \psi}$.

In view of the null form estimate \eqref{NullFormEstimate}, we can reduce \eqref{BilinearBB} to
$$
  I^{\pm}
  \lesssim \norm{\psi}_{{X_+^{s,1/2+\varepsilon}}}
  \norm{\psi'}_{X_\pm^{s,1/2+\varepsilon}},
$$
where
\begin{align*}
  I^{\pm} &= \norm{\int_{\R^{1+2}} \frac{\angles{\xi}^{r-1} \theta_\pm}{\angles{\abs{\tau}-\abs{\xi}}^{1/2-2\varepsilon}} \abs{\widetilde \psi(\lambda,\eta)}
  \abs{\widetilde \psi'(\lambda-\tau,\eta-\xi)} \, d\lambda \, d\eta}_{L^2_{\tau,\xi}},
  \\
  \theta_\pm &= \vangle\bigl(\eta,\pm(\eta-\xi)\bigr).
\end{align*}

Let us right away dispose of the low-frequency case, where $\min(\fixedabs{\eta},\fixedabs{\eta-\xi}) \le 1$ in $I^{\pm}$. Then $\angles{\xi} \sim \angles{\max(\fixedabs{\eta},\fixedabs{\eta-\xi})}$, so assuming $\fixedabs{\eta} \le \fixedabs{\eta-\xi}$, as we may by symmetry, we have
\begin{equation}\label{Ilow}
  I^{\pm}
  \lesssim \norm{\fourierabs{\psi} \cdot \overline{\angles{D}^{r-1} \fourierabs{\psi'}} }_{L^2_{t,x}}
  \le \norm{\fourierabs{\psi}}_{L^\infty_{t,x}} \norm{\angles{D}^{r-1}\fourierabs{\psi'}}_{L^2_{t,x}}
  \lesssim  \norm{\fourierabs{\psi}}_{L^\infty_{t,x}} \norm{\psi'}_{{X_\pm^{s,0}}},
\end{equation}
where we used the assumption $r < 1 + s$. But if the support of $\widetilde \psi(\lambda,\eta)$ is restricted to $\fixedabs{\eta} \le 1$, then $\fixednorm{\fourierabs{\psi}}_{L^\infty_{t,x}} \le C_\varepsilon \fixednorm{\psi}_{{X_\pm^{s,1/2+\varepsilon}}}$ for all $s \in \R$, by the inequalities of Hausdorff-Young and Cauchy-Schwarz.

With the low-frequency case accounted for, we henceforth assume that in $I^{\pm}$,
\begin{equation}\label{FrequencySupport}
  \fixedabs{\eta}, \fixedabs{\eta-\xi} \ge 1.
\end{equation}
It will be convenient to use the notation
\begin{gather*}
  F(\lambda,\eta) = \angles{\eta}^s \angles{\lambda+\fixedabs{\eta}}^{1/2+\varepsilon}
  \abs{\widetilde \psi(\lambda,\eta)},
  \qquad
  G_\pm(\lambda,\eta) = \angles{\eta}^s \angles{\lambda\pm\fixedabs{\eta}}^{1/2+\varepsilon}
  \abs{\widetilde \psi'(\lambda,\eta)},
  \\
  A = \abs{\tau}-\abs{\xi},
  \qquad B = \lambda+\fixedabs{\eta},
  \qquad C_\pm = \lambda-\tau\pm\fixedabs{\eta-\xi},
  \\
  \rho_{+} = \abs{\xi} - \bigabs{ \fixedabs{\eta} - \fixedabs{\eta-\xi}},
  \qquad
  \rho_{-} = \fixedabs{\eta} + \fixedabs{\eta-\xi} - \abs{\xi}.
\end{gather*}
We shall need the easily checked estimates (see \cite{Selberg:2006b})
\begin{equation}\label{ThetaEstimates}
  \theta_{+}^2 \sim \frac{\abs{\xi}
  \rho_{+}}{\fixedabs{\eta}\fixedabs{\eta-\xi}},
  \qquad
  \theta_{-}^2 \sim \frac{(\fixedabs{\eta}+\fixedabs{\eta-\xi})
  \rho_{-}}{\fixedabs{\eta}\fixedabs{\eta-\xi}}
  \sim \frac{\rho_-}{\min(\fixedabs{\eta},\fixedabs{\eta-\xi})}.
\end{equation}
The following estimates will also be needed:
\begin{align}
  \label{rEstimateA}
  \rho_{\pm} &\le 2 \min(\fixedabs{\eta},\fixedabs{\eta-\xi}),
  \\
  \label{rEstimateB}
  \rho_{\pm} &\le \abs{A} + \abs{B} + \abs{C_\pm}.
\end{align}
The first one is the triangle inequality, the second is proved in \cite[Lemma 7]{Selberg:2006b}. 

\subsection{Estimate for $I^+$}
By \eqref{ThetaEstimates}, and recalling \eqref{FrequencySupport}, we have
$$
  I^+ \lesssim
  \norm{\int_{\R^{1+2}} \frac{\angles{\xi}^{r-1/2} \rho_+^{1/2} F(\lambda,\eta) G_+(\lambda-\tau,\eta-\xi)}{\angles{\eta}^{1/2+s} \angles{\eta-\xi}^{1/2+s} \angles{A}^{1/2-2\varepsilon} \angles{B}^{1/2+\varepsilon} \angles{C_+}^{1/2+\varepsilon}} \, d\lambda \, d\eta}_{L^2_{\tau,\xi}}.
$$
By \eqref{rEstimateB} and the fact that $\rho_+ \le \abs{\xi}$, we get $\rho_+^{1/2} \lesssim \abs{A}^{1/2-2\varepsilon} \abs{\xi}^{2\varepsilon} + \abs{B}^{1/2} + \abs{C_+}^{1/2}$, hence
$$
  I^+ \lesssim I^+_1 + I^+_2 + I^+_3,
$$
where
\begin{align*}
  I^+_1 &= \norm{\int_{\R^{1+2}} \frac{\angles{\xi}^{r+2\varepsilon-1/2} F(\lambda,\eta) G_+(\lambda-\tau,\eta-\xi)}{\angles{\eta}^{1/2+s} \angles{\eta-\xi}^{1/2+s} \angles{B}^{1/2+\varepsilon} \angles{C_+}^{1/2+\varepsilon}} \, d\lambda \, d\eta}_{L^2_{\tau,\xi}},
  \\
  I^+_2 &= \norm{\int_{\R^{1+2}} \frac{\angles{\xi}^{r-1/2} F(\lambda,\eta) G_+(\lambda-\tau,\eta-\xi)}{\angles{\eta}^{1/2+s} \angles{\eta-\xi}^{1/2+s} \angles{A}^{1/2-2\varepsilon} \angles{C_+}^{1/2+\varepsilon}} \, d\lambda \, d\eta}_{L^2_{\tau,\xi}},
  \\
  I^+_3 &= \norm{\int_{\R^{1+2}} \frac{\angles{\xi}^{r-1/2} F(\lambda,\eta) G_+(\lambda-\tau,\eta-\xi)}{\angles{\eta}^{1/2+s} \angles{\eta-\xi}^{1/2+s} \angles{A}^{1/2-2\varepsilon} \angles{B}^{1/2+\varepsilon}} \, d\lambda \, d\eta}_{L^2_{\tau,\xi}}.
\end{align*}
In effect, one of the ``hyperbolic'' weights $\angles{A}^{1/2},\angles{B}^{1/2},\angles{C}^{1/2}$ in the denominator in $I^+$ has been traded in for the same amount of ``elliptic'' weights. This happens on account of the null condition, which effectively rules out the bad case where all three hyperbolic weights are $\sim 1$ simultaneously (unless one of the elliptic weights is $\sim 1$).

We only need to estimate $I^+_1$ and $I^+_2$, since $I^+_2$ and $I^+_3$ are symmetrical.

\subsubsection{Estimate for $I^+_1$}\label{Cutoff} We rewrite:
$$
   I^+_1 \simeq \norm{\angles{D}^{r+2\varepsilon-1/2} \left(
   \angles{D}^{-1/2} \fourierabs{\psi}
   \cdot
   \overline{\angles{D}^{-1/2} \fourierabs{\psi'}} \right)  }_{L^2_{t,x}}.
$$
At this point we would like to apply Corollary \ref{BilinearCorollary}, so we need to replace $L^2_{t,x}$ by $L_t^qL_x^2$. This would be trivial if we had compact support in time, but we do not, since $\psi$ has been replaced by $\fourierabs{\psi}$. Nevertheless, there is sufficient decay so that we can apply H\"older's inequality in time. To see this, fix a smooth cutoff function $\chi(t)$ such that $\chi(t) = 1$ for $\abs{t} \le 2$. Denote its Fourier transform by $\widehat \chi(\tau)$, and let $\fourierabs{\chi}$ be the function whose Fourier transform is $\fixedabs{\widehat \chi(\tau)}$. Then $\psi = \chi \psi$, by the support assumption on $\psi$, hence
$$
  \widetilde{\fourierabs{\psi}} = \widetilde{\fourierabs{\chi\psi}}
  \le \widetilde{\fourierabs{\chi}\fourierabs{\psi}}.
$$
Thus, for any $q \ge 2$,
\begin{align*}
   I^+_1 &\le \norm{\angles{D}^{r+2\varepsilon-1/2} \left(
   \angles{D}^{-1/2} ( \fourierabs{\chi}\fourierabs{\psi} )
   \cdot
   \overline{\angles{D}^{-1/2} \fourierabs{\psi'}} \right)  }_{L^2_{t,x}}
   \\
   &= \norm{\fourierabs{\chi} \angles{D}^{r+2\varepsilon-1/2} \left(
   \angles{D}^{-1/2} \fourierabs{\psi}
   \cdot
   \overline{\angles{D}^{-1/2} \fourierabs{\psi'}} \right)  }_{L^2_{t,x}}
   \\
   &\le \norm{\fourierabs{\chi}}_{L_t^p} \norm{\angles{D}^{r+2\varepsilon-1/2} \left(
   \angles{D}^{-1/2} \fourierabs{\psi}
   \cdot
   \overline{\angles{D}^{-1/2} \fourierabs{\psi'}} \right)  }_{L_t^qL_x^2},
\end{align*}
where $1/p + 1/q = 1/2$. Note that $\fixednorm{\fourierabs{\chi}}_{L_t^p} < \infty$, by the Hausdorff-Young inequality. It remains to check, using Corollary \ref{BilinearCorollary}, that
\begin{equation}\label{Iplus1}
   \norm{\angles{D}^{r+2\varepsilon-1/2} \left(
   \angles{D}^{-1/2} \fourierabs{\psi}
   \cdot
   \overline{\angles{D}^{-1/2} \fourierabs{\psi'}} \right)  }_{L_t^qL_x^2}
   \lesssim
   \norm{\psi}_{{X_+^{s,1/2+\varepsilon}}} \norm{\psi'}_{{X_+^{s,1/2+\varepsilon}}},
\end{equation}
for some $q \ge 4$ depending on $r$ and $s$. We divide into two cases: $r < 1/2$ and $r \ge 1/2$.

Assume $r < 1/2$. Then \eqref{Iplus1} follows from Corollary \ref{BilinearCorollary} with $s_1 = s_2 = 1/2+s$, provided we can find $q \ge 4$ such that $s_1+s_2 + 1/2-r-2\varepsilon \ge 1-1/q$, $s_1+s_2 > 1/q$ and $s_1,s_2 <1-1/q$. Written out, the conditions are: 
\begin{align}
  \label{Iplus1A}
  1+2s &\ge \frac{1}{2} - \frac{1}{q} + r + 2\varepsilon,
  \\
  \label{Iplus1B}
  1+2s &> \frac{1}{q}
  \\
  \label{Iplus1C}
  s &< \frac{1}{2} - \frac{1}{q}.
\end{align}
We distinguish three cases: $-1/2 < s \le -3/8$, $-3/8 < s < 1/4$ and $s \ge 1/4$. 

Assume $-1/2 < s \le -3/8$. Then we take $1/q = 1+2s-\varepsilon$, so $4 < q < \infty$, and \eqref{Iplus1B} is satisfied. The condition \eqref{Iplus1C} becomes $s < - 1/6 + \varepsilon/3$, which is satisfied since $s \le -3/8$, and \eqref{Iplus1A} becomes $r+3\varepsilon \le 3/2+4s$, which is satisfied in view of \eqref{rConditions}.

Assume $-3/8 < s < 1/4$. Then we set $q = 4$, so \eqref{Iplus1B} and \eqref{Iplus1C} are certainly satisfied, and \eqref{Iplus1A} reduces to $r+2\varepsilon \le 3/4+2s$, and this is holds by \eqref{rConditions}. 

Finally, assume $s \ge 1/4$. Then there is a lot of room in the estimate. We take $q = 4$ and discard the multiplier $\angles{D}^{r+2\varepsilon-1/2}$ in the left side of \eqref{Iplus1}. Now we take $s_1 = s_2 = 3/8$ in Corollary \ref{BilinearCorollary}. Note that this argument works for $s \ge -1/8$. This concludes the case $r < 1/2$.

Now assume $r \ge 1/2$. This can only happen if $s > -1/8$, in view of \eqref{rConditions}. Note that by symmetry, we may assume $\fixedabs{\eta} \ge \fixedabs{\eta-\xi}$ in $I^+_1$, hence $\angles{\xi} \lesssim \angles{\eta}$. Effectively, we can therefore move the multiplier $\angles{D}^{1/2-r}$ in \eqref{Iplus1} in front of $\fourierabs{\psi}$. Taking $q = 4$, we thus reduce \eqref{Iplus1} to
\begin{equation}\label{Iplus1b}
   \norm{\angles{D}^{r+2\varepsilon-1} \fourierabs{\psi}
   \cdot
   \overline{\angles{D}^{-1/2} \fourierabs{\psi'}} }_{L_t^4L_x^2}
   \lesssim
   \norm{\psi}_{{X_+^{s,1/2+\varepsilon}}} \norm{\psi'}_{{X_+^{s,1/2+\varepsilon}}}.
\end{equation}
To prove this, we apply Corollary \ref{BilinearCorollary}, dividing into two cases: $-1/8 < s < 1/4$ and $s \ge 1/4$. 

Assume $-1/8 < s < 1/4$. Set $s_2=1/2+s$. Then $0 < s_2 < 3/4$, so setting $s_1=3/4-s_2=1/4-s$ and applying Corollary \ref{BilinearCorollary} gives
$$
   \norm{\angles{D}^{r+2\varepsilon-1} \fourierabs{\psi}
   \cdot
   \overline{\angles{D}^{-1/2} \fourierabs{\psi'}} }_{L_t^4L_x^2}
   \lesssim
   \norm{\psi}_{{X_+^{r+2\varepsilon-3/4-s,1/2+\varepsilon}}} \norm{\psi'}_{{X_+^{s,1/2+\varepsilon}}},
$$
which proves \eqref{Iplus1b}, provided $r+2\varepsilon-3/4-s \le s$, i.e., $r+2\varepsilon \le 3/4+2s$, which holds by \eqref{rConditions}.

Now assume $s \ge 1/4$. Set $s_1 = \varepsilon$ and $s_2 = 3/4-\varepsilon$. Then Corollary \ref{BilinearCorollary} gives
$$
   \norm{\angles{D}^{r+2\varepsilon-1} \fourierabs{\psi}
   \cdot
   \overline{\angles{D}^{-1/2} \fourierabs{\psi'}} }_{L_t^4L_x^2}
   \lesssim
   \norm{\psi}_{{X_+^{r+3\varepsilon-1,1/2+\varepsilon}}} \norm{\psi'}_{{X_+^{1/4-\varepsilon,1/2+\varepsilon}}}.
$$
This proves \eqref{Iplus1b}, since $r+3\varepsilon \le 1+s-r$, by \eqref{rConditions}. This concludes the proof of \eqref{Iplus1}.

\subsubsection{Estimate for $I^+_2$}\label{Duality}

Note that if $\abs{A} \ge \abs{\xi}$, then $\rho_+^{1/2} \le \abs{\xi}^{1/2} \le \angles{A}^{1/2-2\varepsilon} \abs{\xi}^{2\varepsilon}$, which means that $I^+$ reduces to $I^+_1$. Therefore, we may assume $\abs{A} \le \abs{\xi}$, which implies
$$
  I^+_2 \lesssim \norm{\int_{\R^{1+2}} \frac{\angles{\xi}^{r+3\varepsilon-1/2} F(\lambda,\eta) G_+(\lambda-\tau,\eta-\xi)}{\angles{\eta}^{1/2+s} \angles{\eta-\xi}^{1/2+s} \angles{A}^{1/2+\varepsilon} \angles{C_+}^{1/2+\varepsilon}} \, d\lambda \, d\eta}_{L^2_{\tau,\xi}}.
$$
We need to prove $I^+_2 \lesssim \norm{F}_{L^2} \norm{G}_{L^2}$, but by duality this is equivalent to the estimate
$$
  \int_{\R^6} \frac{\angles{\xi}^{r+3\varepsilon-1/2} F(\lambda,\eta) G_+(\lambda-\tau,\eta-\xi)H(\tau,\xi)}{\angles{\eta}^{1/2+s} \angles{\eta-\xi}^{1/2+s} \angles{A}^{1/2+\varepsilon} \angles{C_+}^{1/2+\varepsilon}} \, d\lambda \, d\eta \, d\tau \, d\xi
  \lesssim \norm{F}_{L^2} \norm{G_+}_{L^2} \norm{H}_{L^2}
$$
for all $H \in L^2(\R^{1+2})$, $H \ge 0$. For this, it suffices to prove
$$
  \norm{\int_{\R^{1+2}} \frac{\angles{\xi}^{r+3\varepsilon-1/2} G_+(\lambda-\tau,\eta-\xi)H(\tau,\xi)}{\angles{\eta}^{1/2+s} \angles{\eta-\xi}^{1/2+s} \angles{A}^{1/2+\varepsilon} \angles{C_+}^{1/2+\varepsilon}} \, d\tau \, d\xi}_{L^2_{\lambda,\eta}}
  \lesssim \norm{G_+}_{L^2} \norm{H}_{L^2},
$$
which we rewrite as
$$
  \norm{\angles{D}^{-1/2-s} \left(
  \angles{D}^{r+3\varepsilon-1/2} u_\pm
  \cdot
  \angles{D}^{-1/2} \fourierabs{\psi'}
  \right)  }_{L^2_{t,x}}
  \lesssim
   \bignorm{u_\pm}_{{X_\pm^{0,1/2+\varepsilon}}} \norm{\psi'}_{{X_+^{s,1/2+\varepsilon}}},
$$
where
$$
  \widetilde u_\pm(\tau,\xi) = \frac{H_\pm(\tau,\xi)}{\angles{\tau\pm\abs{\xi}}^{1/2+\varepsilon}},
  \quad H_+(\tau,\xi) = \chi_{(-\infty,0)}(\tau)H(\tau,\xi),
  \quad H_-(\tau,\xi) = \chi_{(0,\infty)}(\tau)H(\tau,\xi).
$$
By the trick used for $I^+_1$, we can estimate the $L^2_{t,x}$ norm on the left by the $L_t^qL_x^2$ norm for any $q \ge 2$. Indeed, in the original estimate \eqref{BilinearBB} we can insert the cutoff $\chi$ in front of $\psi$, since $\chi\psi = \psi$, and then we can move $\chi$ in front of $\psi'$. Thus, it suffices to show that for some $q \ge 4$, depending on $s$ and $r$,
\begin{equation}\label{Iplus2}
  \norm{\angles{D}^{-1/2-s} \left(
  \angles{D}^{r+3\varepsilon-1/2} u_\pm
  \cdot
  \angles{D}^{-1/2} \fourierabs{\psi'} \right)  }_{L_t^qL_x^2}
  \lesssim
   \bignorm{u_\pm}_{{X_\pm^{0,1/2+\varepsilon}}} \norm{\psi'}_{{X_+^{s,1/2+\varepsilon}}}.
\end{equation}
We first do the low-frequency case, where $H(\tau,\xi)$ is supported in $\abs{\xi} \ge 1$. Then the left side of \eqref{Iplus2} with $q=\infty$ can be estimated by (we assume $H \ge 0$)
$$
  \norm{\angles{D}^{-1/2-s} \left(
  u_\pm \cdot \angles{D}^{-1/2} \fourierabs{\psi'}
  \right)  }_{L_t^\infty L_x^2}
  \lesssim
  \norm{u_\pm \cdot
  \angles{D}^{-1/2} \fourierabs{\psi'}
  }_{L_t^\infty L_x^2}.
$$
If $-1/2 < s < 1/2$, we apply Corollary \ref{BilinearCorollary} with $s_1 = 1/2+s$, $s_2 = 1/2-s$ and $s_3 = 0$. If $s \ge 1/2$, we can take $s_1 = s_2 = 1/2$ and $s_3 = 0$.

With the low-frequency case out of the way, we assume from now on that $\abs{\xi} \ge 1$ on the support of $H(\tau,\xi)$. Thus, there is effectively no difference between the multipliers $\angles{D}^{r+3\varepsilon-1/2}$ and $\abs{D}^{r+3\varepsilon-1/2}$ when they act on $u_\pm$. In the following we distinguish three cases: $-1/2 < s < 0$, $0 \le s < 1/2$ and $s \ge 1/2$.

First assume $-1/2 < s < 0$. We take $q = 4$, $s_1 = s_3 = 1/2 + s$ and $s_2 = -1/4-2s$. Clearly, $s_1 < 3/4$, and the condition $s_2 < 3/4$ is just $s > -1/2$. The condition $s_1+s_2 > 1/4$ reduces to $s < 0$. Thus Corollary \ref{BilinearCorollary} proves \eqref{Iplus2}, provided $s_2 \le 1/2-r-3/\varepsilon$, i.e., $r+3\varepsilon \le 3/4 + 2s$, which holds by \eqref{rConditions}.

Now assume $0 \le s < 1/2$. Set $1/q = 1/4-s/2-\varepsilon$. Then $4 < q < \infty$, assuming that $2\varepsilon \le 1/2 - s$. Set $s_1 = s_3 = 1/2 + s$ and $s_2 = -1/4-3s/2+\varepsilon$. The condition $s_1 < 1-1/q$ is then the same as $s < 1/2+2\varepsilon$, which is satisfied. The condition $s_2 < 1-1/q$ is nothing else than $s > -1/2$. The condition $s_1+s_2 > 1/q$ reduces to $2\varepsilon > 0$.  Corollary \ref{BilinearCorollary} therefore applies, and \eqref{Iplus2} follows provided $s_2 \le 1/2-r-3/\varepsilon$, i.e., $r+4\varepsilon \le 3/4 + 3s/2 - r$, which holds by \eqref{rConditions}.

Finally, assume $s \ge 1/2$. We separate the cases $r < 1/2$ and $r \ge 1/2$. If $r < 1/2$, we apply Corollary \ref{BilinearCorollary} with $q = \infty$ and, say, $s_1 = 1/2+3\varepsilon$, $s_2 = -3\varepsilon$ and $s_3 = 1/2$. Then \eqref{Iplus2} follows, provided $\varepsilon$ is small enough. Now assume $r \ge 1/2$. Then using the triangle inequality we can estimate the left side of \eqref{Iplus2} by a sum of two terms, namely
\begin{align*}
  J_1 &= \norm{\angles{D}^{-1-s+r+3\varepsilon} \left(
  u_\pm \cdot \angles{D}^{-1/2} \fourierabs{\psi'}
  \right) }_{L_t^qL_x^2},
  \\
  J_2 &= \norm{\angles{D}^{-1/2-s} \left( u_\pm 
  \cdot \angles{D}^{-1+r+3\varepsilon} \fourierabs{\psi'}
  \right) }_{L_t^qL_x^2}.
\end{align*}
Now apply Corollary \ref{BilinearCorollary} with $q = \infty$. Choose $4\varepsilon \le 1+s-r$, in accordance with \eqref{rConditions}. Then in $J_1$ we can take $s_1 = 1-\varepsilon$, $s_2 = 0$ and $s_3 = \varepsilon$, and in $J_2$ we can take $s_1 = \varepsilon$, $s_2 = 0$ and $s_3 = 1-\varepsilon$. This concludes the proof of \eqref{Iplus2}.

\subsection{Estimate for $I^-$}\label{MinusReduction}

Note first that if $\fixedabs{\eta} \ll \fixedabs{\eta-\xi}$, then $\abs{\xi} \sim \fixedabs{\eta-\xi}$, so by \eqref{ThetaEstimates},
$$
  \theta_-^2 \sim \frac{\abs{\xi}\rho_-}{\fixedabs{\eta}\fixedabs{\eta-\xi}}.
$$
The same is true if $\fixedabs{\eta} \gg \fixedabs{\eta-\xi}$ or $\abs{\xi} \sim \fixedabs{\eta} \sim \fixedabs{\eta-\xi}$. So in all these cases, we have the same estimate for $\theta_-$ as for $\theta_+$, and moreover we have $\rho_-^{1/2} \lesssim \angles{A}^{1/2-2\varepsilon} \abs{\xi}^{2\varepsilon} + \angles{B}^{1/2} + \angles{C_-}^{1/2}$, by \eqref{rEstimateA} and \eqref{rEstimateB}, so the analysis of $I^+$ in the previous subsection then applies also to $I^-$. Thus, it suffices to consider $I^-$ in the case where
\begin{equation}\label{FrequencySupport2}
  \abs{\xi} \ll \fixedabs{\eta} \sim \fixedabs{\eta-\xi},
\end{equation}
which we assume from now on. Recalling also \eqref{FrequencySupport}, we then have, by \eqref{ThetaEstimates},
$$
  I^- \lesssim
  \norm{\int_{\R^{1+2}} \frac{\angles{\xi}^{r-1} \rho_-^{1/2} F(\lambda,\eta) G_-(\lambda-\tau,\eta-\xi)}{\angles{\eta}^{1/4+s} \angles{\eta-\xi}^{1/4+s} \angles{A}^{1/2-2\varepsilon} \angles{B}^{1/2+\varepsilon} \angles{C_-}^{1/2+\varepsilon}} \, d\lambda \, d\eta}_{L^2_{\tau,\xi}}.
$$

Note that for $I^+$ we were able to get an estimate for $s > -1/2$, and $r$ satisfying \eqref{rConditions}, but now we have to assume $s > -1/4$, since the weights $\angles{\eta}$ and $\angles{\eta-\xi}$ are taken to the power $1/4+s$ instead of $1/2+s$. The first condition in \eqref{rConditions} then becomes redundant.

By \eqref{rEstimateA} and \eqref{rEstimateB}, we have $\rho_-^{1/2} \lesssim \abs{A}^{1/2-2\varepsilon} \fixedabs{\eta}^{\varepsilon} \fixedabs{\eta-\xi}^{\varepsilon} + \abs{B}^{1/2} + \abs{C_-}^{1/2}$, so
$$
  I^- \lesssim I^-_1 + I^-_2 + I^-_3,
$$
where
\begin{align*}
  I^-_1 &= \norm{\int_{\R^{1+2}} \frac{\angles{\xi}^{r-1} F(\lambda,\eta) G_-(\lambda-\tau,\eta-\xi)}{\angles{\eta}^{1/4+s-\varepsilon} \angles{\eta-\xi}^{1/4+s-\varepsilon} \angles{B}^{1/2+\varepsilon} \angles{C_-}^{1/2+\varepsilon}} \, d\lambda \, d\eta}_{L^2_{\tau,\xi}},
  \\
  I^-_2 &= \norm{\int_{\R^{1+2}} \frac{\angles{\xi}^{r-1} F(\lambda,\eta) G_-(\lambda-\tau,\eta-\xi)}{\angles{\eta}^{1/4+s} \angles{\eta-\xi}^{1/4+s} \angles{A}^{1/2-2\varepsilon} \angles{C_-}^{1/2+\varepsilon}} \, d\lambda \, d\eta}_{L^2_{\tau,\xi}},
  \\
  I^-_3 &= \norm{\int_{\R^{1+2}} \frac{\angles{\xi}^{r-1} F(\lambda,\eta) G_-(\lambda-\tau,\eta-\xi)}{\angles{\eta}^{1/4+s} \angles{\eta-\xi}^{1/4+s} \angles{A}^{1/2-2\varepsilon} \angles{B}^{1/2+\varepsilon}} \, d\lambda \, d\eta}_{L^2_{\tau,\xi}}.
\end{align*}

\subsubsection{Estimate for $I^-_1$}\label{DifficultCase}

Proceeding as in the estimate for $I^+_1$, we can reduce to proving
\begin{equation}\label{Iminus1}
   \norm{\angles{D}^{r-1} \left(
   \angles{D}^{\varepsilon-1/4} \fourierabs{\psi}
   \cdot
   \overline{\angles{D}^{\varepsilon-1/4} \fourierabs{\psi'}} \right) }_{L_t^qL_x^2}
   \lesssim
   \norm{\psi}_{{X_+^{s,1/2+\varepsilon}}} \norm{\psi'}_{{X_-^{s,1/2+\varepsilon}}},
\end{equation}
for some $q \ge 4$. However, this estimate fails for $s < -1/8$. To see this, take $s_1 = s_2 = 1/4+s-\varepsilon$ and consider the conditions in Corollary \ref{BilinearCorollary}, which as remarked are sharp up to endpoints. The condition $1-r+s_1+s_2 \ge 1-1/q$ becomes $2\varepsilon \le 1/2+1/q+2s-r$, forcing $q=4$, since $r$ can be arbitrarily close to $3/4+2s$, by \eqref{rConditions}. The condition $s_1+s_2 \ge 1/q$ then forces $s \ge -1/8 + \varepsilon$.

Fortunately, there is a way around this difficulty. The solution is to keep the $L^2$ norm, instead of passing to a higher $L^q$ norm in time. Recall \eqref{FrequencySupport2}, which says that we are in the high-high frequency case with output at low frequency. Moreover, the interaction is of the type $(+,+)$, because of the conjugation in the second factor. Indeed, observe that $\angles{C_-} = \angles{\lambda-\tau-\fixedabs{\eta-\xi}} = \angles{\tau-\lambda+\fixedabs{\xi-\eta}}$, so we can rewrite $I^-_1$ in a form which is manifestly of type $(+,+)$, namely
$$
  I^-_1 = \norm{\int_{\R^{1+2}} \frac{\angles{\xi}^{r-1} F(\lambda,\eta) G'(\tau-\lambda,\xi-\eta)}{\angles{\eta}^{1/4+s-\varepsilon} \angles{\eta-\xi}^{1/4+s-\varepsilon} \angles{\lambda+\fixedabs{\eta}}^{1/2+\varepsilon} \angles{\tau-\lambda+\fixedabs{\xi-\eta}}^{1/2+\varepsilon}} \, d\lambda \, d\eta}_{L^2_{\tau,\xi}},
$$
where we have set $G'(\lambda,\eta) = G_-(-\lambda,-\eta)$. We now see from Corollary \ref{BilinearCorollary2}, using \eqref{FrequencySupport2}, that
$$
  I^-_1 \lesssim \bignorm{F}_{L^2} \norm{G'}_{L^2}
  = \norm{\psi}_{{X_+^{s,1/2+\varepsilon}}} \norm{\psi'}_{{X_-^{s,1/2+\varepsilon}}},
$$
provided $2\varepsilon \le 1+2s-r$, $s \le \varepsilon + 3/8$ and $0 < s+1/4-\varepsilon$. These conditions are satisfied for $-1/4 < s \le 3/8$, in view of \eqref{rConditions}. If $s > 3/8$, we can estimate $\angles{\xi}^{r-1} \le \angles{\xi}^s$, since $r<1+s$, and apply Corollary \ref{BilinearCorollary} with $q = 4$ (in fact, this works for $s > 1/4$). This concludes the proof for $I^-_1$.

\subsubsection{Estimate for $I^-_2$}

We may assume $\abs{A} \ll \fixedabs{\eta}+\fixedabs{\eta-\xi}$, since otherwise $I^- \lesssim I^-_1$, by \eqref{rEstimateA}. Recalling also \eqref{FrequencySupport2}, we then get
$$
  I^-_2 \lesssim \norm{\int_{\R^{1+2}} \frac{\angles{\xi}^{r-1} F(\lambda,\eta) G_-(\lambda-\tau,\eta-\xi)}{\angles{\eta-\xi}^{1/2+2s-3\varepsilon} \angles{A}^{1/2+\varepsilon} \angles{C_-}^{1/2+\varepsilon}} \, d\lambda \, d\eta}_{L^2_{\tau,\xi}}.
$$
Estimating by duality as in the proof of $I^+_2$, and using the cutoff argument, we reduce to proving
$$
  \norm{\angles{D}^{r-1} u_\pm
  \cdot
  \angles{D}^{3\varepsilon-1/2-s} \fourierabs{\psi'}  }_{L_t^4L_x^2}
  \lesssim
  \bignorm{u_\pm}_{{X_\pm^{0,1/2+\varepsilon}}} \norm{\psi'}_{{X_-^{s,1/2+\varepsilon}}}.
$$
Assuming $r < 1$, this follows from Corollary \ref{BilinearCorollary}, provided $1-r+1/2+2s-3\varepsilon \ge 3/4$, i.e., $3\varepsilon \le 3/4+2s-r$, which is compatible with \eqref{rConditions}. On the other hand, if $r \ge 1$, then using \eqref{FrequencySupport2}, we see that
$$
  I^-_2 \lesssim \norm{\int_{\R^{1+2}} \frac{F(\lambda,\eta) G_-(\lambda-\tau,\eta-\xi)}{\angles{\xi}^{a/2}\angles{\eta-\xi}^{a/2} \angles{A}^{1/2+\varepsilon} \angles{C_-}^{1/2+\varepsilon}} \, d\lambda \, d\eta}_{L^2_{\tau,\xi}},
$$
where $a = 3/2+2s-r-3\varepsilon$. Hence we reduce to, by duality and cutoff,
$$
  \norm{\angles{D}^{-a/2} u_\pm
  \cdot
  \angles{D}^{-a/2+s} \fourierabs{\psi'}  }_{L_t^4L_x^2}
  \lesssim
  \bignorm{u_\pm}_{{X_\pm^{0,1/2+\varepsilon}}} \norm{\psi'}_{{X_-^{s,1/2+\varepsilon}}},
$$
and this follows from Corollary \ref{BilinearCorollary}, provided $a \ge 3/4$, i.e., $3\varepsilon \le 3/4+2s-r$.

\subsubsection{Estimate for $I^-_3$} The argument used for $I^-_2$ applies also here. In fact, by a change of variables, $I^-_3$ can be transformed to $I^-_2$, except that the minus sign changes to a plus sign: we will have $C_+$ instead of $C_-$. However, the argument used for $I^-_2$ is not sensitive to a change of sign. This concludes the proof of Theorem \ref{Thm2}.

\section{Proof of Theorem \ref{Thm3}}\label{Thm3Proof}

We take $\pm_1 = +$ and write $\pm_2 = \pm$. We assume $s \in \R$ and
\begin{equation}\label{rConditions2}
  r > \max\left(-s,\frac{1}{4}-\frac{s}{2},\frac{1}{4}+\frac{s}{2},s\right).
\end{equation}
Let $\varepsilon > 0$. As remarked after the statement of Theorem \ref{Thm3}, it suffices to prove the estimate \eqref{BilinearAAA}. We assume that $\psi,\psi' \in \mathcal S(\R^{1+2})$ and $\psi$ is supported in $[-2,2] \times \R^2$.

Using \eqref{NullFormEstimate}, we reduce \eqref{BilinearAAA} to
\begin{equation}\label{IIplusminus}
  I^{\pm}
  \lesssim \norm{\psi}_{{X_+^{s,1/2+\varepsilon}}}
  \norm{\psi'}_{X_\pm^{-s,1/2-2\varepsilon}},
\end{equation}
where now
$$
  I^{\pm} = \norm{\int_{\R^{1+2}} \frac{\theta_\pm}{\angles{\xi}^{r} \angles{\abs{\tau}-\abs{\xi}}^{1/2+\varepsilon}} \abs{\widetilde \psi(\lambda,\eta)}
  \abs{\widetilde \psi'(\lambda-\tau,\eta-\xi)} \, d\lambda \, d\eta}_{L^2_{\tau,\xi}},
$$
and $\theta_\pm = \vangle\bigl(\eta,\pm(\eta-\xi)\bigr)$ as before. We proceed as in the proof of Theorem \ref{Thm2}, and use the notation introduced there, except that now
$$
  G_\pm(\lambda,\eta) = \angles{\eta}^{-s} \angles{\lambda\pm\fixedabs{\eta}}^{1/2-2\varepsilon}
  \abs{\widetilde \psi'(\lambda,\eta)}.
$$

Consider the low-frequency case, $\min(\fixedabs{\eta},\fixedabs{\eta-\xi}) \le 1$. Then $\angles{\xi} \sim \angles{\max(\fixedabs{\eta},\fixedabs{\eta-\xi})}$. If $\fixedabs{\eta} \le 1$, then \eqref{Ilow} applies, with $r-1$ replaced by $-r$, and $s$ replaced by $-s$, hence we need $-r \le -s$, which holds by \eqref{rConditions2}. On the other hand, if $\fixedabs{\eta-\xi} \le 1$, then using $-r \le s$, which holds by \eqref{rConditions2}, and using also the Sobolev type estimate
$$
  \norm{u}_{L_t^\infty L_x^2} \simeq \norm{\widehat u}_{L_t^\infty L_\xi^2} \le \norm{\widehat u}_{L_\xi^2 L_t^\infty} \le \norm{\widetilde u}_{L_\xi^2 L_\tau^1} \le C_\varepsilon \norm{u}_{X_\pm^{0,1/2+\varepsilon}},
$$
we can write
$$
  I^{\pm}
  \lesssim \norm{\angles{D}^{-r} \fourierabs{\psi} \cdot \overline{\fourierabs{\psi'}} }_{L^2_{t,x}}
  \le \norm{\angles{D}^{-r} \fourierabs{\psi}}_{L_t^\infty L_x^2} 
   \norm{\fourierabs{\psi'}}_{L_t^2 L_x^\infty}
   \lesssim \norm{\psi}_{{X_\pm^{s,1/2+\varepsilon}}} \norm{\fourierabs{\psi'}}_{L_t^2 L_x^\infty}.
$$
But since we are assuming now that $\widetilde \psi'(\lambda-\tau,\eta-\xi)$ vanishes unless $\fixedabs{\eta-\xi} \le 1$, we have $\fixednorm{\fourierabs{\psi'}}_{L_t^2L_x^\infty} \lesssim \fixednorm{\psi'}_{{X_\pm^{-s,0}}}$ for all $s \in \R$, by Sobolev embedding. This concludes the proof of the low-frequency case.

From now on we assume the high-frequency case, so that in $I^{\pm}$,
\begin{equation}\label{HighFrequency}
  \fixedabs{\eta}, \fixedabs{\eta-\xi} \ge 1.
\end{equation}

\subsection{Estimate for $I^+$}
From \eqref{ThetaEstimates} and \eqref{HighFrequency} we get
$$
  I^+ \lesssim
  \norm{\int_{\R^{1+2}} \frac{\rho_+^{1/2} F(\lambda,\eta) G_+(\lambda-\tau,\eta-\xi)}{\angles{\xi}^{r-1/2} \angles{\eta}^{1/2+s} \angles{\eta-\xi}^{1/2-s} \angles{A}^{1/2+\varepsilon} \angles{B}^{1/2+\varepsilon} \angles{C_+}^{1/2-2\varepsilon}} \, d\lambda \, d\eta}_{L^2_{\tau,\xi}}.
$$
By \eqref{rEstimateA} and \eqref{rEstimateB}, $\rho_+^{1/2} \lesssim \abs{A}^{1/2} + \abs{B}^{1/2} + \abs{C_+}^{1/2-2\varepsilon} \fixedabs{\eta-\xi}^{2\varepsilon}$, hence
$
  I^+ \lesssim I^+_1 + I^+_2 + I^+_3,
$
where
\begin{align*}
  I^+_1 &= \norm{\int_{\R^{1+2}} \frac{F(\lambda,\eta) G_+(\lambda-\tau,\eta-\xi)}{\angles{\xi}^{r-1/2} \angles{\eta}^{1/2+s} \angles{\eta-\xi}^{1/2-s} \angles{B}^{1/2+\varepsilon} \angles{C_+}^{1/2-2\varepsilon}} \, d\lambda \, d\eta}_{L^2_{\tau,\xi}},
  \\
  I^+_2 &= \norm{\int_{\R^{1+2}} \frac{F(\lambda,\eta) G_+(\lambda-\tau,\eta-\xi)}{\angles{\xi}^{r-1/2} \angles{\eta}^{1/2+s} \angles{\eta-\xi}^{1/2-s} \angles{A}^{1/2+\varepsilon} \angles{C_+}^{1/2-2\varepsilon}} \, d\lambda \, d\eta}_{L^2_{\tau,\xi}},
  \\
  I^+_3 &= \norm{\int_{\R^{1+2}} \frac{F(\lambda,\eta) G_+(\lambda-\tau,\eta-\xi)}{\angles{\xi}^{r-1/2} \angles{\eta}^{1/2+s} \angles{\eta-\xi}^{1/2-s-2\varepsilon} \angles{A}^{1/2+\varepsilon} \angles{B}^{1/2+\varepsilon}} \, d\lambda \, d\eta}_{L^2_{\tau,\xi}}.
\end{align*}

\subsubsection{Estimate for $I^+_1$} Note that the sum of the exponents of the ``elliptic'' weights $\angles{\xi}$, $\angles{\eta}$ and $\angles{\eta-\xi}$ in $I^+_1$ is $1/2+r$. If $r > 1/2$, it therefore follows by a Sobolev type estimate, see \cite[Proposition A.1]{Selberg:2002b}, that $I^+_1 \lesssim \norm{F}_{L^2} \norm{G_+}_{L^2}$. Henceforth we may therefore assume
$r \le 1/2$.
Thus, we may restrict to
$-1/2 < s < 1/2$,
since $\abs{s} \ge 1/2$ forces $r > 1/2$, by \eqref{rConditions2}.

Next, note that if $\abs{C_+} \ge \fixedabs{\eta-\xi}$, then the sum of the elliptic exponents increases to $1+r-2\varepsilon$, which exceeds $1$, since we assume $r > 1/4$. Since we still have the weight $\angles{B}$ to a power greater than $1/2$, it now follows that $I^+_1 \lesssim \norm{F}_{L^2} \norm{G_+}_{L^2}$, by the Sobolev estimate referred to above.

In view of the above reductions, we can from now on assume
\begin{equation}\label{IIplus1A}
  r \le 1/2,
  \qquad -1/2 < s < 1/2,
  \qquad \abs{C_+} \le \fixedabs{\eta-\xi}.
\end{equation}
Thus, $\angles{\xi}^{1/2-r}  \lesssim \angles{\eta}^{1/2-r} + \angles{\eta-\xi}^{1/2-r}$, and we may assume
$$
  \norm{\psi'}_{{X_+^{-s-3\varepsilon,1/2+\varepsilon}}}
  \le \norm{\psi'}_{{X_+^{-s,1/2-2\varepsilon}}}.
$$
By the same cutoff argument as in subsection \ref{Cutoff}, we then reduce to proving
\begin{align}
  \label{IIplus1B}
  \norm{\angles{D}^{-r} \fourierabs{\psi}
  \cdot
  \overline{\angles{D}^{-1/2} \fourierabs{\psi'}} }_{L_t^4L_x^2}
  &\lesssim
  \norm{\psi}_{{X_+^{s,1/2+\varepsilon}}} \norm{\psi'}_{{X_+^{-s-3\varepsilon,1/2+\varepsilon}}},
  \\
  \label{IIplus1C}
  \norm{\angles{D}^{-1/2} \fourierabs{\psi}
  \cdot
  \overline{\angles{D}^{-r} \fourierabs{\psi'}} }_{L_t^4L_x^2}
  &\lesssim
  \norm{\psi}_{{X_+^{s,1/2+\varepsilon}}} \norm{\psi'}_{{X_+^{-s-3\varepsilon,1/2+\varepsilon}}}.
\end{align}
Now we apply Corollary \ref{BilinearCorollary} with $q = 4$ and $s_3 = 0$. Sufficient conditions for this to work are, for \eqref{IIplus1B} and \eqref{IIplus1C} respectively,
\begin{alignat*}{3}
  r+s &> 0,& \qquad 1/2-s-3\varepsilon &> 0,& \qquad r+1/2-3\varepsilon &\ge 3/4,
  \\
  1/2+s &> 0,& \qquad r-s-3\varepsilon &> 0,& \qquad r+1/2-3\varepsilon &\ge 3/4.
\end{alignat*}
All these conditions are satisfied, in view of \eqref{rConditions2} and the assumption $-1/2 < s < 1/2$. This concludes the proof for $I^+_1$.

\subsubsection{Estimate for $I^+_2$} By the same reductions as in the previous subsection, we may assume \eqref{IIplus1A}, and by the duality and cutoff argument used in subsection \ref{Duality}, we then reduce to proving
\begin{align}
  \label{IIplus2A}
  \norm{\angles{D}^{-r-s} \left(
  u_\pm
  \cdot
  \angles{D}^{-1/2} \fourierabs{\psi'}
  \right)  }_{L_t^q L_x^2}
  &\lesssim
  \bignorm{u_\pm}_{{X_\pm^{0,1/2+\varepsilon}}} \norm{\psi'}_{{X_+^{-s-3\varepsilon,1/2+\varepsilon}}},
  \\
  \label{IIplus2B}
  \norm{\angles{D}^{-1/2-s} \left(
  u_\pm
  \cdot
  \angles{D}^{-r} \fourierabs{\psi'}
  \right)  }_{L_t^q L_x^2}
  &\lesssim
  \bignorm{u_\pm}_{{X_\pm^{0,1/2+\varepsilon}}} \norm{\psi'}_{{X_+^{-s-3\varepsilon,1/2+\varepsilon}}},
\end{align}
for some $q \ge 4$.

For \eqref{IIplus2A}, we apply Corollary \ref{BilinearCorollary} with $s_1 = 0$ and $s_2 = 1/2-s-3\varepsilon$. Thus, we need
\begin{gather}
  \label{IIplus2A1}
  r+\frac{1}{2}-3\varepsilon \ge 1-\frac{1}{q},
  \\
  \label{IIplus2A2}
  \frac{1}{q} < \frac{1}{2}-s-3\varepsilon < 1-\frac{1}{q}.
\end{gather}
We divide into three cases: $-1/4 < s < 1/4$, $-1/2 < s \le -1/4$ and $1/4 \le s < 1/2$.

If $-1/4 < s < 1/4$, then \eqref{IIplus2A2} is satisfied with $q=4$, and then \eqref{IIplus2A1} reduces to $3\varepsilon \le r-1/4$, which is in accordance with \eqref{rConditions2}. If $-1/2 < s \le -1/4$, then the left inequality in \eqref{IIplus2A2} is certainly satisfied, and to optimize we choose $q \ge 4$ as small as possible, but such that the right inequality in \eqref{IIplus2A2} is satisfied. Thus, we set $1/q=s+1/2$, so $4 \le q < \infty$. Then \eqref{IIplus2A2} holds, and \eqref{IIplus2A1} reduces to $3\varepsilon \le r+s$, which is compatible with \eqref{rConditions2}. Finally, if $1/4 \le s < 1/2$, we set $1/q=-s+1/2$. Then \eqref{IIplus2A2} again holds, and \eqref{IIplus2A1} becomes $3\varepsilon \le r-s$, in agreement with \eqref{rConditions2}.

For \eqref{IIplus2B}, we apply Corollary \ref{BilinearCorollary} with $s_1 = 0$ and $s_2 = r-s-3\varepsilon$. We then require \eqref{IIplus2A1} to hold, and also
\begin{equation}
  \label{IIplus2B2}
  \frac{1}{q} < r-s-3\varepsilon < 1-\frac{1}{q}.
\end{equation}
We distinguish the same cases as above.

If $-1/4 \le s < 0$, set $q = 4$. Then \eqref{IIplus2A1} holds, and \eqref{IIplus2B2} reduces to $-1/4-3\varepsilon < s < 1/4-3\varepsilon$, which is also satisfied. If $-1/2 < s \le -1/4$, we take $1/q=s+1/2$. Again, \eqref{IIplus2A1} is satisfied, and \eqref{IIplus2B2} becomes $1/2+2s+3\varepsilon < r < 1/2+3\varepsilon$, which holds since $1/4 < r \le 1/2$ by \eqref{rConditions2} and \eqref{IIplus1A}. Finally, assume $1/4 \le s < 1/2$. Then the right inequality in \eqref{IIplus2B2} is certainly satisfied, so we are left with \eqref{IIplus2A1}  and the left inequality in \eqref{IIplus2B2}, which we can restate as two lower bounds for $r$:
$$
  r \ge \frac{1}{2}-\frac{1}{q} + 3\varepsilon,
  \qquad
  r >s+\frac{1}{q} + 3\varepsilon.
$$
To optimize, we take $1/q=1/4-s/2$. Then we get the condition $r > 1/4 + s/2 + 3\varepsilon$, which is compatible with \eqref{rConditions2}. This concludes the proof for $I^+_2$.

\subsubsection{Estimate for $I^+_3$} Again, we may assume \eqref{IIplus1A}. Thus, $-1/2 < s < 1/2$, so replacing $s$ by $-s$, and noting also the invariance of \eqref{rConditions2}, we then reduce $I^+_3$ to $I^+_2$. In particular, we need here the condition $r > 1/4 - s/2$. 

\subsection{Estimate for $I^-$}

By the argument in subsection \ref{MinusReduction}, we may assume
$$
\abs{\xi} \ll \fixedabs{\eta} \sim \fixedabs{\eta-\xi}
$$
in $I^-$. Combining this with \eqref{ThetaEstimates} and \eqref{HighFrequency}, we get
$$
  I^- \lesssim
  \norm{\int_{\R^{1+2}} \frac{\rho_-^{1/2} F(\lambda,\eta) G_-(\lambda-\tau,\eta-\xi)}{\angles{\xi}^{r} \angles{\eta}^{1/4} \angles{\eta-\xi}^{1/4} \angles{A}^{1/2+\varepsilon} \angles{B}^{1/2+\varepsilon} \angles{C_+}^{1/2-2\varepsilon}} \, d\lambda \, d\eta}_{L^2_{\tau,\xi}}.
$$
By \eqref{rEstimateA} and \eqref{rEstimateB}, $\rho_-^{1/2} \lesssim \abs{A}^{1/2} + \abs{B}^{1/2} + \abs{C_-}^{1/2-2\varepsilon} \fixedabs{\eta-\xi}^{2\varepsilon}$, hence
$
  I^- \lesssim I^-_1 + I^-_2 + I^-_3,
$
where
\begin{align*}
  I^-_1 &= \norm{\int_{\R^{1+2}} \frac{F(\lambda,\eta) G_-(\lambda-\tau,\eta-\xi)}{\angles{\xi}^{r} \angles{\eta}^{1/4} \angles{\eta-\xi}^{1/4} \angles{B}^{1/2+\varepsilon} \angles{C_-}^{1/2-2\varepsilon}} \, d\lambda \, d\eta}_{L^2_{\tau,\xi}},
  \\
  I^-_2 &= \norm{\int_{\R^{1+2}} \frac{F(\lambda,\eta) G_-(\lambda-\tau,\eta-\xi)}{\angles{\xi}^{r} \angles{\eta}^{1/4} \angles{\eta-\xi}^{1/4} \angles{A}^{1/2+\varepsilon} \angles{C_-}^{1/2-2\varepsilon}} \, d\lambda \, d\eta}_{L^2_{\tau,\xi}},
  \\
  I^-_3 &= \norm{\int_{\R^{1+2}} \frac{F(\lambda,\eta) G_-(\lambda-\tau,\eta-\xi)}{\angles{\xi}^{r} \angles{\eta}^{1/4-\varepsilon} \angles{\eta-\xi}^{1/4-\varepsilon} \angles{A}^{1/2+\varepsilon} \angles{B}^{1/2+\varepsilon}} \, d\lambda \, d\eta}_{L^2_{\tau,\xi}},
\end{align*}
so $s$ no longer appears. Using the cutoff argument, and possibly duality, we reduce to Corollary \ref{BilinearCorollary} with $q = 4$, and we conclude that \eqref{IIplusminus} holds with the minus sign, provided $r > 1/4$, which is satisfied on account of \eqref{rConditions2}. This concludes the proof of Theorem \ref{Thm3}.

\section{Optimality of Theorems \ref{Thm2} and \ref{Thm3}}\label{Optimality}

Here we prove the optimality, up to endpoint cases, of the conditions on $s$ and $r$ in Theorems \ref{Thm2} and \ref{Thm3}. We first consider Theorem \ref{Thm2}, commenting on Theorem \ref{Thm3} at the end of this section. Thus, we first prove that \eqref{BilinearBB} fails unless
\begin{align}
  \label{sConditionsNonsharp}
  s &\ge -\frac{1}{4},
  \\
  \label{rConditionsNonsharp}
  r &\le \min\left(\frac{3}{4}+2s,\frac{3}{4}+\frac{3s}{2},1+s\right).
\end{align}
Since we ignore endpoints, we may set $\varepsilon = 0$ in \eqref{BilinearBB}. The following counterexamples all depend on a parameter $L \gg 1$ going to infinity. Note that
\begin{multline}\label{Product1}
  \norm{\innerprod{\beta\Pi_{+}(D) \psi}{\Pi_{\pm}(D)\psi'}}_{H^{r-1,-1/2}}
  \\
  =
  \norm{\int_{\R^{1+2}} \frac{\innerprod{\Pi_{\pm}(\eta-\xi)\beta\Pi_{+}(\eta)\widetilde \psi(\lambda,\eta)}{\widetilde \psi'(\lambda-\tau,\eta-\xi)}}{\angles{\xi}^{1-r}\angles{\abs{\tau}-\abs{\xi}}^{1/2-2\varepsilon}} \, d\lambda \, d\eta}_{L^2_{\tau,\xi}}.
\end{multline}
In each counterexample we choose either the plus or the minus sign, and we choose $A,B,C \subset \R^2$, depending on $L$ and concentrated along the $\xi_1$-direction, with the property
\begin{equation}\label{ABCproperty}
  \eta \in A, \,\, \xi \in C \implies \eta-\xi \in B.
\end{equation}
We then construct $\psi$ and $\psi'$ depending on $L$, one of which is supported in $[-2,2] \times \R^2$, such that
\begin{equation}\label{LowerBound}
  \frac{\norm{\innerprod{\beta\Pi_{+}(D) \psi}{\Pi_{\pm}(D)\psi'}}_{H^{r-1,-1/2}}}
  {\norm{\psi}_{X_+^{s,1/2}}\norm{\psi'}_{X_+^{s,1/2}}}
  \ge \frac{1}{CL^{\delta(r,s)}},
\end{equation}
where $C$ is independent of $L$, and $\delta(r,s)$ depends on the choice of $A,B,C$. In each case we get a necessary condition $\delta(r,s) \ge 0$. We shall use the notation $\mathbb{1}_{(\cdot)}$ for the function whose value is $1$ if the condition $(\cdot)$ in the subscript is satisfied, and $0$ otherwise. 

\subsection{Necessity of \eqref{rConditionsNonsharp}.}\label{rOptimality}

For this, we take the plus sign in \eqref{Product1}, and choose $\widetilde\psi$ and $\widetilde\psi'$ to be characteristic functions of slabs cut out of a thickened null hyperplane $\tau+\xi_1 = O(1)$, multiplied by an eigenvector of $\Pi_+(\xi)$ (see \eqref{Projection}),
\begin{equation}\label{vPlus}
  v_+(\xi) = \bigl[ 1, \hat\xi_1+i\hat\xi_2\bigr]^T,
  \qquad \text{where} \quad \hat\xi \equiv \frac{\xi}{\abs{\xi}}.
\end{equation}
Assuming $A,B,C$ have been chosen, we set
\begin{align}
  \label{Psi1}
  \widetilde \psi(\lambda,\eta) &= \mathbb{1}_{\lambda+\eta_1 = O(1)} \mathbb{1}_{\eta \in A} v_+(\eta),
  \\\label{Psi2}
  \widetilde \psi'(\lambda-\tau,\eta-\xi) &= \mathbb{1}_{\lambda-\tau+\eta_1-\xi_1 = O(1)} \mathbb{1}_{\eta -\xi \in B} v_+(\eta-\xi),
\end{align}
and we restrict the $L^2$ norm to the region
\begin{equation}\label{Output1}
  \tau+\xi_1 = O(1), \quad \xi \in C.
\end{equation}
The fact that we work in a thickened hyperplane trivializes the convolution structure, since obviously,
\begin{equation}\label{Obvious}
  \tau+\xi_1=O(1), \,\,\lambda+\eta_1=O(1) \implies \lambda-\tau+\eta_1-\xi_1=O(1).
\end{equation}
Using this and \eqref{ABCproperty}, we get from \eqref{Product1}, restricting the $L^2$ norm to the region \eqref{Output1},
\begin{equation}\label{Product2}
\begin{aligned}
  &\norm{\innerprod{\beta\Pi_{+}(D) \psi}{\Pi_+(D)\psi'}}_{H^{r-1,-1/2}}
  \\
  &\qquad=
  \norm{\int_{\R^{1+2}} \frac{\innerprod{\beta v_{+}(\eta)}{v_+(\eta-\xi)}}{\angles{\xi}^{1-r}\angles{\abs{\tau}-\abs{\xi}}^{1/2}}
  \mathbb{1}_{\left\{
  \scriptstyle\eta \in A
  \atop
  \scriptstyle\lambda+\eta_1 = O(1)
  \right\}
  }
  \mathbb{1}_{\left\{
    \scriptstyle\eta-\xi \in B
    \atop
  \scriptstyle\lambda-\tau+\eta_1-\xi_1 = O(1)
  \right\}
  }
  \, d\lambda \, d\eta}_{L^2_{\tau,\xi}}
  \\
  &\qquad\ge
  \norm{\int_{\R^{1+2}} \frac{\innerprod{\beta v_{+}(\eta)}{v_+(\eta-\xi)}}{\angles{\xi}^{1-r}\angles{\abs{\tau}-\abs{\xi}}^{1/2}}
  \mathbb{1}_{\left\{
  \scriptstyle\eta \in A
  \atop
  \scriptstyle\lambda+\eta_1 = O(1)
  \right\}
  }
  \mathbb{1}_{\left\{
  \scriptstyle\xi \in C
  \atop
  \scriptstyle\tau+\xi_1 = O(1)
  \right\}
  }
  \, d\lambda \, d\eta}_{L^2_{\tau,\xi}}
  \\
  &\qquad\ge
  \norm{\int_{\R^{1+2}} \frac{\im\innerprod{\beta v_{+}(\eta)}{v_+(\eta-\xi)}}{\angles{\xi}^{1-r}\angles{\abs{\tau}-\abs{\xi}}^{1/2}}
  \mathbb{1}_{\left\{
  \scriptstyle\eta \in A
  \atop
  \scriptstyle\lambda+\eta_1 = O(1)
  \right\}
  }
  \mathbb{1}_{\left\{
  \scriptstyle\xi \in C
  \atop
  \scriptstyle\tau+\xi_1 = O(1)
  \right\}
  }
  \, d\lambda \, d\eta}_{L^2_{\tau,\xi}}.
\end{aligned}
\end{equation}
But
\begin{equation}\label{Eigenproduct}
  \innerprod{\beta v_{+}(\eta)}{v_+(\zeta)}
  = 1 - \hat\eta\cdot\hat\zeta + i\hat\eta\wedge\hat\zeta,
  \qquad \text{where}
  \quad\hat\eta\wedge\hat\zeta=\hat\eta_1\hat\zeta_2 - \hat\eta_2\hat\zeta_1,
\end{equation}
hence
\begin{equation}\label{Imagpart}
  \im\innerprod{\beta v_{+}(\eta)}{v_+(\eta-\xi)} = \pm\sin\theta_+ \sim \pm\theta_+,
  \qquad
  \text{where} \quad \theta_+ = \vangle(\eta,\eta-\xi).
\end{equation}
The choice of sign in front of $\sin\theta_+$ depends on the orientation: The sign is $+$ if rotating $\eta$ counterclockwise through the angle $\theta_+$ makes it line up with $\eta-\xi$. But the sets $A,B,C$ will be chosen so that the orientation of the pair $(\eta,\eta-\xi)$ is fixed. From \eqref{Product2} we therefore conclude:
\begin{multline}\label{Product3}
  \norm{\innerprod{\beta\Pi_{+}(D) \psi}{\Pi_+(D)\psi'}}_{H^{r-1,-1/2}} \ge I^+,
  \\
  \text{where}
  \quad
  I^+
  =
  \norm{\int_{\R^{1+2}} \frac{\theta_+}{\angles{\xi}^{1-r}\angles{\abs{\tau}-\abs{\xi}}^{1/2}}
  \mathbb{1}_{\left\{
  \scriptstyle\eta \in A
  \atop
  \scriptstyle\lambda+\eta_1 = O(1)
  \right\}
  }
  \mathbb{1}_{\left\{
  \scriptstyle\xi \in C
  \atop
  \scriptstyle\tau+\xi_1 = O(1)
  \right\}
  }
  \, d\lambda \, d\eta}_{L^2_{\tau,\xi}}.
\end{multline}

So far we ignored the requirement that $\psi$ or $\psi'$ be supported in $[-2,2] \times \R^2$, but this is easily fixed. Let $\chi(t)$ be a smooth, even cutoff function such that $0 \le \chi \le 1$, $\chi(t) = 1$ for $\abs{t} \le 1$ and $\chi(t) = 0$ for $\abs{t} \ge 2$. Since $\chi$ is even, its Fourier transform $\widehat \chi(\tau)$ is real-valued. Moreover, $\widehat \chi(\tau) \sim 1$ for $\abs{\tau} = O(1)$, if we interpret $O(1)$ to mean that $\abs{\tau} \le \delta$ for some constant $\delta > 0$ determined by $\chi$. For convenience we choose to localize $\psi'$ rather than $\psi$. Note that
\begin{equation}\label{Psi2cutoff}
  \widetilde{\chi\psi'}(\lambda-\tau,\eta-\xi)
  =  f(\lambda-\tau+\eta_1-\xi_1) \mathbb{1}_{\eta -\xi \in B} v_+(\eta-\xi),
\end{equation}
where
$$
  f(\tau) = \int_{\R} \widehat \chi(\mu)\mathbb{1}_{\tau-\mu = O(1)} \, d\mu.
$$
So the only difference from \eqref{Psi2} is that $\mathbb{1}_{\lambda-\tau+\eta_1-\xi_1}$ has been replaced by $f(\lambda-\tau+\eta_1-\xi_1)$. But from \eqref{Obvious} and the fact that $f(\tau) \sim 1$ for $\tau = O(1)$, we then conclude that \eqref{Product3} still holds if we replace $\psi'$ by $\chi\psi'$. On the other hand, by \eqref{SimpleEstimate} we have
$$
  \norm{\psi'}_{X_+^{s,1/2}} \gtrsim \norm{\chi\psi'}_{X_+^{s,1/2}},
$$
and we conclude that if \eqref{LowerBound} holds without the cutoff, then it also holds with the cutoff. In what follows we can therefore ignore the cutoff function.

We now construct the counterexamples, by choosing the sets $A,B,C$. Note that in $I^+$,
\begin{equation}\label{Slabs}
\begin{gathered}
  \eta \in A, \qquad \xi \in C, \qquad \eta-\xi \in B,
  \\
  \lambda + \eta_1 = O(1),
  \qquad \tau + \xi_1 = O(1),
  \qquad \lambda-\tau + \eta_1-\xi_1 = O(1).
\end{gathered}
\end{equation}

\subsubsection{Necessity of $r \le 3/4+2s$}

\begin{figure}[h]
   \centering
   \includegraphics{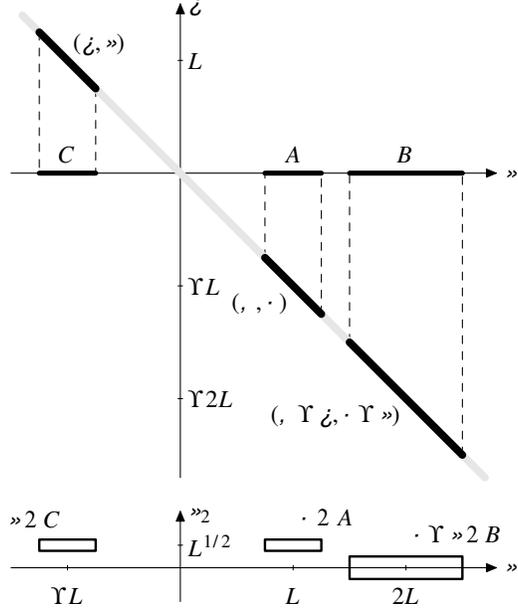}
   \caption{Geometry of the sets used to prove necessity of $r \le 3/4+2s$.}
   \label{fig:2}
\end{figure}

Here we consider high frequencies interacting to give output at high frequency. Set (see Figure \ref{fig:2})
\begin{align*}
  A &= \left\{ \xi \in \R^2 : \abs{\xi_1 - L} \le L/4, \,\, \abs{\xi_2 - L^{1/2}} \le L^{1/2}/4 \right\},
  \\
  B &= \left\{ \xi \in \R^2 : \abs{\xi_1 - 2L} \le L/2, \,\, \abs{\xi_2} \le L^{1/2}/2 \right\},
  \\
  C &= \left\{ \xi \in \R^2 : \abs{\xi_1 + L} \le L/4, \,\, \abs{\xi_2 - L^{1/2}} \le L^{1/2}/4 \right\}.
\end{align*}

Then \eqref{ABCproperty} holds. Recalling \eqref{Slabs}, we see that
$$
  \theta_+ = \vangle(\eta,\eta-\xi) \sim \frac{1}{L^{1/2}},
  \qquad \abs{\xi}, \fixedabs{\eta}, \fixedabs{\eta-\xi} \sim L,
$$
and
\begin{equation}\label{NearCone1}
  \lambda+\fixedabs{\eta} = \lambda+\eta_1+\fixedabs{\eta}-\eta_1
  = \lambda+\eta_1 + \frac{\eta_2^2}{\fixedabs{\eta}+\eta_1} = O(1).
\end{equation}
Similarly,
\begin{equation}\label{NearCone2}
  \lambda-\tau+\fixedabs{\eta-\xi} = O(1),
  \qquad
  \abs{\tau}-\abs{\xi} = \tau-\abs{\xi} = O(1).
\end{equation}
Thus, denoting by $\abs{A}$ the area of $A$, we see that
$$
  I^+ \sim \frac{\abs{A}\abs{C}^{1/2}}{L^{1/2+1-r}}
  \qquad \text{and} \qquad
  \norm{\psi}_{X_+^{s,1/2}}, \norm{\psi'}_{X_+^{s,1/2}}
  \sim L^s \abs{A}^{1/2}.
$$
Since $\abs{A} = \abs{C} \sim L^{3/2}$, we conclude that \eqref{LowerBound} holds with
$\delta(r,s) = 3/4-r+2s$, proving the necessity of $r \le 3/4+2s$.

\subsubsection{Necessity of $r \le 3/4+3s/2$}\label{HLH1}

\begin{figure}[h]
   \centering
   \includegraphics{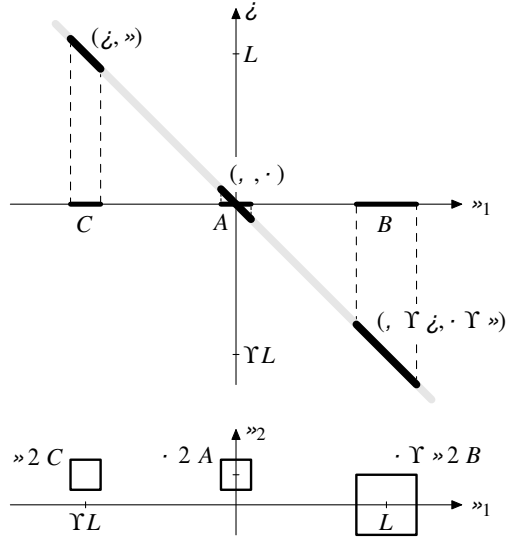}
   \caption{Geometry of sets used to prove necessity of $r \le 3/4+3s/2$ and $r \le 1+s$, the squares $A,B,C$ having sides $\sim L^{1/2}$ or $\sim 1$, respectively.}
   \label{fig:3}
\end{figure}

Here we use a high/low frequency interaction with output at high frequency. Define (see Figure \ref{fig:3})
\begin{align*}
  A &= \left\{ \xi \in \R^2 : \abs{\xi_1} \le L^{1/2}/2, \,\, \abs{\xi_2 - 1} \le L^{1/2}/2 \right\},
  \\
  B &= \left\{ \xi \in \R^2 : \abs{\xi_1 - L} \le L^{1/2}, \,\, \abs{\xi_2} \le L^{1/2} \right\},
  \\
  C &= \left\{ \xi \in \R^2 : \abs{\xi_1 + L} \le L^{1/2}/2, \,\, \abs{\xi_2 - 1} \le L^{1/2}/2 \right\}.
\end{align*}

Then
$\theta_+ = \vangle(\eta,\eta-\xi) \sim 1$,
$\fixedabs{\eta} \sim L^{1/2}$
and $\abs{\xi}, \fixedabs{\eta-\xi} \sim L$.
Further, \eqref{NearCone2} still holds, whereas the calculation in \eqref{NearCone1} shows that
$\lambda+\fixedabs{\eta} \sim L^{1/2}$,
since $\fixedabs{\eta}+\eta_1 \ge \eta_2 - \eta_1 \ge L^{1/2}/2$. Thus,
$$
  I^+ \sim \frac{\abs{A}\abs{C}^{1/2}}{L^{1-r}},
  \qquad
  \norm{\psi}_{X_+^{s,1/2}} \sim L^{s/2+1/4} \abs{A}^{1/2},
  \qquad \norm{\psi'}_{X_+^{s,1/2}}
  \sim L^s \abs{B}^{1/2}.
$$
But $\abs{A}, \abs{B}, \abs{C} \sim L$, hence \eqref{LowerBound} holds with
$\delta(r,s) = 3/4-r+3s/2$, proving the necessity of $r \le 3/4+3s/2$.

\subsubsection{Necessity of $r \le 1+s$}\label{HLH2}

The configuration is the same as in the previous subsection, except that the squares $A,B,C$ now have side length $\sim 1$. We set
\begin{align*}
  A &= \left\{ \xi \in \R^2 : \abs{\xi_1} \le 1/2, \,\, \abs{\xi_2 - 1} \le 1/2 \right\},
  \\
  B &= \left\{ \xi \in \R^2 : \abs{\xi_1 - L} \le 1, \,\, \abs{\xi_2} \le 1 \right\},
  \\
  C &= \left\{ \xi \in \R^2 : \abs{\xi_1 + L} \le 1/2, \,\, \abs{\xi_2 - 1} \le 1/2 \right\}.
\end{align*}
Then $\theta_+ \sim 1$, $\fixedabs{\eta} \sim 1$, $\abs{\xi}, \fixedabs{\eta-\xi} \sim L$,
and \eqref{NearCone1} holds, since $\fixedabs{\eta}+\eta_1 \ge \eta_2 - \eta_1 \ge 1/2$. Since \eqref{NearCone2} also holds, we conclude:
$$
  I^+ \sim \frac{\abs{A}\abs{C}^{1/2}}{L^{1-r}},
  \qquad
  \norm{\psi}_{X_+^{s,1/2}} \sim \abs{A}^{1/2},
  \qquad \norm{\psi'}_{X_+^{s,1/2}}
  \sim L^s \abs{B}^{1/2}.
$$
But $\abs{A}, \abs{B}, \abs{C} \sim 1$, so \eqref{LowerBound} holds with
$\delta(r,s) = 1-r+s$, proving necessity of $r \le 1+s$.

\subsection{Necessity of \eqref{sConditionsNonsharp}.}

\begin{figure}[h]
   \centering
   \includegraphics{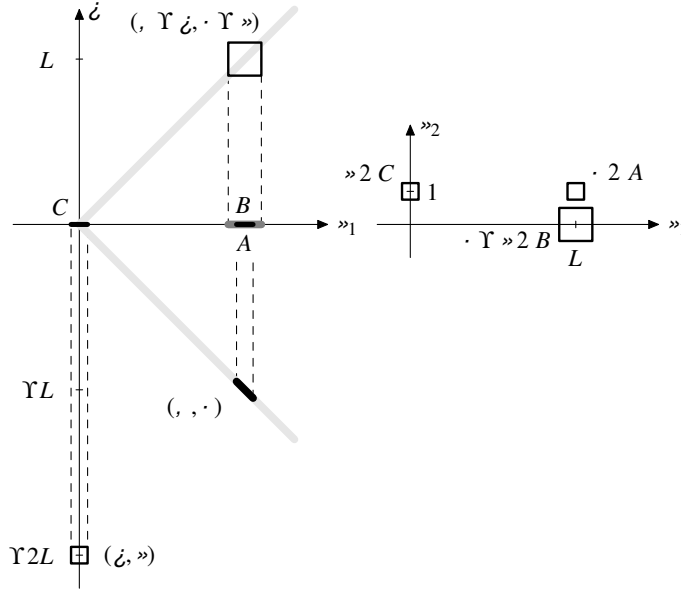}
   \caption{Geometry of sets used to prove necessity of $s\ge-1/4$. The squares $A,B,C$ have side length $\sim 1$.}
   \label{fig:4}
\end{figure}

Here we consider high frequencies interacting to give output at low frequency, and we choose the minus sign in \eqref{Product1}. Set (see Figure \ref{fig:4})
\begin{align*}
  A &= \left\{ \xi \in \R^2 : \abs{\xi_1 - L} \le 1/4, \,\, \abs{\xi_2 - 1} \le 1/4 \right\},
  \\
  B &= \left\{ \xi \in \R^2 : \abs{\xi_1 - L} \le 1/2, \,\, \abs{\xi_2} \le 1/2 \right\},
  \\
  C &= \left\{ \xi \in \R^2 : \abs{\xi_1} \le 1/4, \,\, \abs{\xi_2 - 1} \le 1/4 \right\}.
\end{align*}

In \eqref{Product1} we now restrict the integration to
$$
  \eta \in A, \qquad \lambda + \fixedabs{\eta} = O(1),
  \qquad
  \xi \in C, \qquad \tau + 2L = O(1),
$$
which implies
$$
  \eta-\xi \in B, \qquad \lambda-\tau - \fixedabs{\eta-\xi} = \lambda+\fixedabs{\eta}-\tau-2L+L-\fixedabs{\eta}+L-\fixedabs{\eta-\xi}= O(1),
$$
since $L-\fixedabs{\eta} = L-\eta_1-\eta_2^2/(\fixedabs{\eta}+\eta_1) = O(1)$ and, similarly, $L-\fixedabs{\eta-\xi} = O(1)$. Now set
\begin{align*}
  \widetilde \psi(\lambda,\eta) &= \mathbb{1}_{\lambda+\fixedabs{\eta} = O(1)} \mathbb{1}_{\eta \in A} v_+(\eta),
  \\
  \widetilde \psi'(\lambda-\tau,\eta-\xi) &= \mathbb{1}_{\lambda-\tau-\fixedabs{\eta-\xi} = O(1)} \mathbb{1}_{\eta -\xi \in B} v_-(\eta-\xi),
\end{align*}
where $v_-(\xi) = v_+(-\xi)$ and $v_+(\xi)$ is given by \eqref{vPlus}. Thus, $v_-(\xi)$ is an eigenvector of $\Pi_-(\xi) = \Pi_+(-\xi)$. Since
$
  \theta_- = \vangle(\eta,\xi-\eta) \sim 1,
$
we then get, arguing as in \eqref{Product2}, and using \eqref{Eigenproduct},
\begin{multline*}
  \norm{\innerprod{\beta\Pi_{+}(D) \psi}{\Pi_-(D)\psi'}}_{H^{r-1,-1/2}} \ge I^-,
  \\
  \text{where}
  \quad
  I^-
  =
  \norm{\int_{\R^{1+2}} \frac{1}{\angles{\xi}^{1-r}\angles{\abs{\tau}-\abs{\xi}}^{1/2}}
  \mathbb{1}_{\left\{
  \scriptstyle\eta \in A
  \atop
  \scriptstyle\lambda+\fixedabs{\eta} = O(1)
  \right\}
  }
  \mathbb{1}_{\left\{
  \scriptstyle\xi \in C
  \atop
  \scriptstyle\tau+2L = O(1)
  \right\}
  }
  \, d\lambda \, d\eta}_{L^2_{\tau,\xi}}.
\end{multline*}
Since $\abs{\xi} \sim 1$, $\fixedabs{\eta}, \fixedabs{\eta-\xi} \sim L$ and $\abs{\tau}-\abs{\xi} \sim \abs{\tau} \sim L$, we see that
$$
  I^- \sim \frac{\abs{A} \abs{C}^{1/2}}{L^{1/2}},
  \qquad
  \norm{\psi}_{X_+^{s,1/2}} \sim L^{s} \abs{A}^{1/2},
  \qquad \norm{\psi'}_{X_+^{s,1/2}}
  \sim L^s \abs{B}^{1/2}.
$$
But $\abs{A}, \abs{B}, \abs{C} \sim 1$, hence \eqref{LowerBound} holds with
$\delta(r,s) = 1/2+2s$, proving necessity of $s \ge -1/4$.
By the same reasoning as in Section \ref{rOptimality}, we can also include a time cutoff (then we must rescale the squares $A,B,C$ to have side length $\sim \delta$, where $\delta > 0$ is a small number depending on the cutoff).

\subsection{Optimality in Theorem \ref{Thm3}}

Here we prove that \eqref{BilinearAAA} fails unless
$$
  r \ge \min\left(-s,\frac{1}{4}-\frac{s}{2},\frac{1}{4}+\frac{s}{2},s\right).
$$
In fact, the counterexample in subsection \ref{HLH1} shows that $r \ge 1/4-s/2$ is necessary, hence, by symmetry, also $r \ge 1/4+s/2$, and the example in \ref{HLH2} shows the necessity of $r \ge s$, hence also $r \ge -s$.

\section{Proof of Theorem \ref{Thm4}}\label{Thm4Proof}

We use the following sharp dyadic decomposition in frequency space. By a \emph{dyadic number} we mean a number of the form $2^j$, where $j \in \Z$. We let $\kappa$, $\lambda$ and $\mu$ be dyadic numbers. Given $f \in \mathcal S(\R^2)$, we denote by $f_\lambda$ the function whose Fourier transform is $\chi_{A_{\lambda}} \widehat f$, where $A_{\lambda}$ is the annulus $\lambda \le  \abs{\xi} \le 2\lambda$. We also use a decomposition of $\R^2$ into almost disjoint dyadic squares $Q = [j\mu,(j+1)\mu] \times [k\mu,(k+1)\mu]$ of side length $\mu$, where $j,k \in \Z$; we shall refer to these squares as \emph{$\mu$-squares}. We denote by $f_\lambda^Q$ the function whose Fourier transform is $\chi_{A_{\lambda} \cap Q} \widehat f$.
Recall that $S_\pm(t) = e^{\mp i t \abs{D}}$ is the free propagator for $-i\partial_t \pm \abs{D}$. We let $f,g \in \mathcal S(\R^2)$ and write, for any combination of signs,
$$
  u(t) = S_\pm(t)f,
  \qquad
  v(t) = S_\pm(t)g.
$$
We also write $u_\lambda(t) = S_\pm(t) f_\lambda$ and $u_\lambda^Q(t) = S_\pm(t) f_\lambda^Q$, and similarly for $v$.

An exponent pair $(q,r)$ is said to be \emph{wave admissible} if $2/q + 1/r \le 1/2$ and $r \neq \infty$. For $(q,r)$ wave admissible, the Strichartz estimate holds:
\begin{equation}\label{LinearStrichartz}
  \bignorm{u}_{L_t^q L_x^r} \le C_{q,r}\norm{f}_{\dot H^s} \qquad \left( s = 1 - \frac{2}{r} - \frac{1}{q} \right).
\end{equation}
Note that if $r = 4$, then all $q \ge 8$ are admissible, and the following bilinear version holds:
\begin{equation}\label{BilinearStrichartzA}
  \bignorm{uv}_{L_t^q L_x^2} \le C_{q,\delta} \norm{f}_{\dot H^\delta}  \norm{g}_{\dot H^{1-1/q-\delta}}
  \qquad (4 \le q \le \infty, \quad 0 < \delta < 1-1/q).
\end{equation}
This follows from \eqref{LinearStrichartz} using H\"older's inequality with suitable wave admissible pairs $(q_1,r_1)$ and $(q_2,r_2)$. In fact, by interpolation it suffices to consider $q = 4$ and $q = \infty$.

As proved in \cite{Klainerman:1999}, the constant in \eqref{LinearStrichartz} can be improved if $\widehat f$ is supported in a small set: If $\widehat f$ is supported in a $\mu$-square at distance $\sim \lambda$ from the origin, then the constant can be replaced by $C_{q,r} (\mu/\lambda)^{1/2-1/r}$. In fact, we shall need only the case $(q,r) = (8,4)$:
\begin{equation}\label{ImprovedStrichartzA}
  \bignorm{u_\lambda^Q}_{L_t^8 L_x^4} \le C \left( \frac{\mu}{\lambda} \right)^{1/4} \bignorm{f_\lambda^Q}_{\dot H^{3/8}}
  \sim \mu^{1/4} \lambda^{1/8} \bignorm{f_\lambda^Q}_{L^2},
\end{equation}
where $Q$ is a $\mu$-square.

We also need the analogous small-support improvement of the Sobolev inequality $\norm{f}_{L_x^4} \le C \norm{f}_{\dot H^{1/2}}$. By the Hausdorff-Young inequality,
$$
  \bignorm{f_\lambda^Q}_{L_x^4} \le \norm{\chi_{A_{\lambda} \cap Q} \widehat f}_{L_\xi^{4/3}}
  \le \norm{\chi_{A_{\lambda} \cap Q}}_{L_\xi^{4}} \norm{\chi_{A_{\lambda} \cap Q}\widehat f}_{L^2_\xi}
  \le \mu^{1/2} \bignorm{f_\lambda^Q}_{L^2_x},
$$
where $Q$ is a $\mu$-square. Hence, $\bignorm{S_\pm(t) f_\lambda^Q}_{L_t^\infty L_x^4} \le \mu^{1/2} \bignorm{S_\pm(t) f_\lambda^Q}_{L_t^\infty L_x^2} = \mu^{1/2} \bignorm{f_\lambda^Q}_{L^2_x}$, and interpolation between this and \eqref{ImprovedStrichartzA} gives
\begin{equation}\label{ImprovedStrichartzB}
  \bignorm{u_\lambda^Q}_{L_t^q L_x^4} \le
  C \mu^{1/2-2/q} \lambda^{1/q} \bignorm{f_\lambda^Q}_{L^2}
  \qquad (8 \le q \le \infty).
\end{equation}
Using this and \eqref{BilinearStrichartzA} we now prove Theorem \ref{Thm4}.

Write $u = \sum_{\kappa} u_{\kappa}$, $v = \sum_{\lambda} v_{\lambda}$, $u_{\ll \lambda} = \sum_{\kappa \le \lambda/4} u_{\lambda_1}$ and $v_{\ll \kappa} = \sum_{\lambda \le \kappa/4} v_{\lambda}$, where $\kappa$ and $\lambda$ are dyadic numbers. Then 
$$
  uv = \sum_{\lambda} u_{\ll \lambda} v_{\lambda}
  + \sum_{\kappa} \sum_{\kappa/2 \le \lambda \le 2\kappa} u_{\kappa} v_{\lambda}
  + \sum_{\kappa} u_{\kappa} v_{\ll \kappa}
  = S_1 + S_2 + S_3.
$$
The estimate for $S_2$ reduces to \eqref{ImprovedStrichartzB}, by the argument used in the proof of Theorem 4 in Appendix A of \cite{Klainerman:1999}, so we restrict our attention to $S_1$ and $S_3$. By symmetry, it suffices to consider $S_1$.

Since the Fourier transform of $u_{\ll \lambda} v_{\lambda}$ is supported in the annulus $\lambda/2 \le \abs{\xi} \le 4\lambda$, we have, by orthogonality,
$$
  \norm{\abs{D}^{-s_3} S_1}_{L^2_x}
  \sim \Bigl( \sum_{\lambda} \norm{\abs{D}^{-s_3} (u_{\ll \lambda} v_{\lambda})}_{L^2_x}^2 \Bigr)^{1/2}
  \sim \Bigl( \sum_{\lambda} \lambda^{-2s_3} \bignorm{u_{\ll \lambda} v_{\lambda}}_{L^2_x}^2 \Bigr)^{1/2},
$$
so by Minkowski's integral inequality,
$$
  \bignorm{\abs{D}^{-s_3} S_1}_{L_t^q L_x^2}
  \lesssim \Bigl( \sum_{\lambda} \lambda^{-2s_3} \bignorm{u_{\ll \lambda} v_{\lambda}}_{L_t^q L_x^2}^2 \Bigr)^{1/2}.
$$
Now we apply \eqref{BilinearStrichartzA}, considering separately the cases $s_1 > 0$ and $s_1 \le 0$. If $s_1 > 0$, we get, noting that $s_2+s_3=1-1/q-s_1 > 0$,
$$
  \bignorm{\abs{D}^{-s_3} S_1}_{L_t^q L_x^2}
  \lesssim \Bigl( \sum_{\lambda} \lambda^{-2s_3} \bignorm{f_{\ll \lambda}}_{\dot H^{s_1}}^2
  \norm{g_{\lambda}}_{\dot H^{s_2+s_3}}^2 \Bigr)^{1/2}
  \lesssim  \bignorm{f}_{\dot H^{s_1}}
  \Bigl( \sum_{\lambda} \norm{g_{\lambda}}_{\dot H^{s_2}}^2 \Bigr)^{1/2}.
$$
If $s_1 \le 0$, then $s_2 > 1/q - s_1 \ge 0$, and since also $s_1 + s_3 > 0$, we get
\begin{align*}
  \bignorm{\abs{D}^{-s_3} S_1}_{L_t^q L_x^2}
  &\lesssim \Bigl( \sum_{\lambda} \lambda^{-2s_3} \bignorm{f_{\ll \lambda}}_{\dot H^{s_1+s_3}}^2
  \norm{g_{\lambda}}_{\dot H^{s_2}}^2 \Bigr)^{1/2}
  \\
  &\sim \Bigl( \sum_{\lambda} \sum_{\kappa \le \lambda/4} \left(\frac{\kappa}{\lambda}\right)^{2s_3} \bignorm{f_{\kappa}}_{\dot H^{s_1}}^2
  \norm{g_{\lambda}}_{\dot H^{s_2}}^2 \Bigr)^{1/2}
  \lesssim  \bignorm{f}_{\dot H^{s_1}}
  \Bigl( \sum_{\lambda} \norm{g_{\lambda}}_{\dot H^{s_2}}^2 \Bigr)^{1/2}.
\end{align*}
To get the last inequality we used the fact that $s_3 > 0$, since $s_3 = (1-1/q-s_2)-s_1$.

\end{document}